\documentclass[10pt, journal]{IEEEtran}
\pdfoutput=1

\makeatletter
\usepackage[top= 1 in, bottom=0.75 in, left=0.75 in, right=0.75 in]{geometry}

\usepackage{cite}

\usepackage{graphics,color}
\usepackage{graphicx}
\usepackage{epsfig, psfrag}
\usepackage{subfigure}
\usepackage{epstopdf}
\usepackage{psfrag}
\usepackage[normalem]{ulem}

\usepackage{amsmath}
\usepackage{latexsym}
\usepackage{amssymb}
\usepackage{amsthm} 

\usepackage{placeins}

\usepackage{url}
\usepackage{xspace}
\usepackage{soul}

\newtheorem{theorem}{Theorem}
\newtheorem{lemma}{Lemma}

\newtheorem{corollary}{Corollary}
\newtheorem{assumption}{Assumption}
\newtheorem{fact}{Fact}
\theoremstyle{definition}
\newtheorem{definition}{Definition}
\newtheorem{remarks}{Remark}
\newtheorem{example}{Example}

\def\expe{\mathbb{E}}   

\def\limsup{\mathop{\rm limsup}}

\def\P{\mathsf{P}}

\def\mc{\mathcal}
\def\mbf{\mathbf}

\newfont{\boldlarge}{msbm10 scaled 1100}

\newcommand\Ber{\mathsf{Ber}}

\newcommand{\kl}[2]{D\left( \left. #1 \right\| #2 \right)}

\newcommand\defeq{\stackrel{\triangle}{=}}
\newcommand{\dist}[3]{f_{#1}\left(#2; #3\right)}
\newcommand{\samp}[2]{X_{#1}^{(#2)}}
\newcommand{\msg}[2]{Y_{#1}^{(#2)}}
\newcommand{\bel}[2]{b_{#1}^{(#2)}}
\newcommand{\est}[2]{q_{#1}^{(#2)}}

\newcommand{\ignore}[1]{}
\newcommand{\estc}[1]{\mathbf{\tilde{q}_{#1}^{(t)}}}

\newcommand{\add}[1]{{\color{black}{#1}}}
\newcommand{\remove}[1]{}

\renewcommand{\qed}{\nobreak \ifvmode \relax \else
      \ifdim\lastskip<1.5em \hskip-\lastskip
      \hskip1.5em plus0em minus0.5em \fi \nobreak
      \vrule height0.2em width0.5em depth0.4em\fi}

\begin{document}

\title{Social Learning and Distributed Hypothesis Testing}

\author{Anusha~Lalitha, 
        Tara~Javidi,~\IEEEmembership{Senior~Member,~IEEE,}
        and~Anand~Sarwate,~\IEEEmembership{Member,~IEEE}
\thanks{This paper was presented in part in~\cite{DHT_ISIT2014, LalithaJ:15allerton, DHT_asilomar2015}.}
\thanks{A. Lalitha and T. Javidi are with the Department of Electrical and Computer Engineering, University of California San Diego, La Jolla, CA 92093, USA. (e-mail: alalitha@ucsd.edu; tjavidi@ucsd.edu). } 
\thanks{A. Sarwate is with the Department of Electrical and Computer Engineering, Rutgers, The State University of New Jersey, 94 Brett Road, Piscataway, NJ 08854 , USA. (e-mail: asarwate@ece.rutgers.edu).}
}

\maketitle

\begin{abstract}
This paper considers a problem of distributed hypothesis testing and social learning.  Individual nodes in a network receive noisy local (private) observations whose distribution is parameterized by a discrete parameter (hypotheses).  The marginals of the joint observation distribution conditioned on each hypothesis are known locally at the nodes, but the true parameter/hypothesis is not known.  
An update rule is analyzed in which nodes first perform a  Bayesian update of their belief (distribution estimate) of each hypothesis based on their local observations, communicate these updates to their neighbors, and then perform a ``non-Bayesian'' linear consensus using the log-beliefs of their neighbors.  Under mild assumptions, we show that the belief of any node on a wrong hypothesis converges to zero exponentially fast, and the exponential rate of learning is characterized by the nodes' influence of the network and the divergences between the observations' distributions. For a broad class of observation statistics which includes distributions with unbounded support such as Gaussian mixtures, we show that rate of rejection of wrong hypothesis satisfies a large deviation principle i.e., the probability of sample paths on which the rate of rejection of wrong hypothesis deviates from the mean rate vanishes exponentially fast and we characterize the rate function in terms of the nodes' influence of the network and the local observation models.
\end{abstract}


\section{Introduction}

Learning in a distributed setting is more than a phenomenon of social networks; it is also an engineering challenge for 
networked system designers.  For instance, in today's data networks, many 
applications need estimates of certain parameters:
file-sharing systems need to know the distribution
of (unique) documents shared by their users, internet-scale information retrieval
systems need to deduce the criticality 
of various data items, and monitoring networks need to
compute aggregates in a duplicate-insensitive manner.  
Finding scalable, efficient, and accurate methods for
computing such metrics (e.g. number of documents
in the network, sizes of database relations, distributions of
data values) is of critical value in a wide array of network applications.

\begin{figure}
  \centering
    \includegraphics[width=0.4\textwidth]{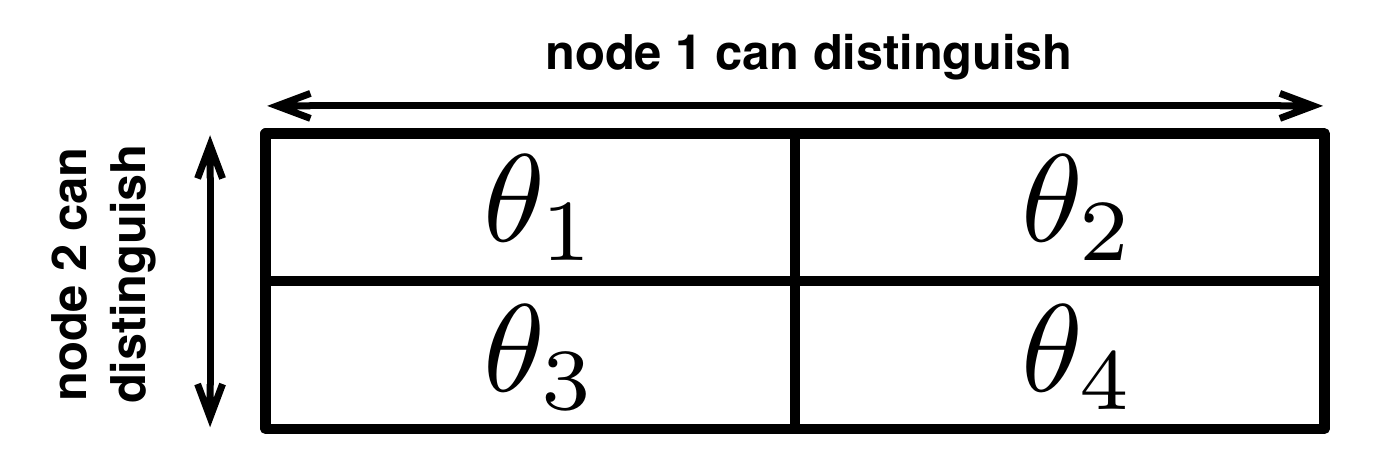}
  \caption{Example of a parameter space in which no node can identify the true parameter.  There are 4 parameters, $\{\theta_1, \theta_2, \theta_3, \theta_4\}$, and 2 nodes.  The node 1 has $\dist{1}{\cdot}{\theta_1} = \dist{1}{\cdot}{\theta_3}$ and $\dist{1}{\cdot}{\theta_2} = \dist{1}{\cdot}{\theta_4}$, and the node 2 has $\dist{2}{\cdot}{\theta_1} = \dist{2}{\cdot}{\theta_2}$ and $\dist{2}{\cdot}{\theta_3} = \dist{2}{\cdot}{\theta_4}$.}
\label{fig:grid}
\end{figure}

We consider a network of individuals sample local observations (over time) governed by an unknown true hypothesis $\theta^*$ taking values in a finite discrete set $\Theta$. We model the $i$-th node's distribution (or local channel, or likelihood function) of the observations conditioned on the true hypothesis by $\dist{i}{\cdot}{\theta^*}$ from a collection $\{ \dist{i}{\cdot}{\theta} : \theta \in \Theta \}$. Nodes neither have access to each others' observation nor the joint distribution of observations across all nodes in the network. A simple two-node example is illustrated in Figure \ref{fig:grid} -- one node can only learn the column in which the true hypothesis lies, and the other can only learn the row.  The local observations of a given node are not sufficient to recover the underlying hypothesis in isolation. \add{In this paper we study a learning rule that enables the nodes to learn the unknown true hypothesis based on message passing between one hop neighbors (local communication) in the network. In particular, each node performs a local Bayesian update and send its belief vectors (message) to its neighbors. After receiving the messages from the neighbors each node performs a consensus averaging on a reweighting of the \textit{log beliefs}. Our result shows that under our learning rule each node can reject the wrong hypothesis exponentially fast.}

We show that the rate of rejection of wrong hypothesis is the weighted sum of Kullback-Leibler (KL) divergences between likelihood function of the true parameter and the likelihood function of the wrong hypothesis, where the sum is over the nodes in the network and the weights are the nodes' influences as dictated by the learning rule. 
\add{
Furthermore, we show that the probability of sample paths on which the rate of rejection deviates from the mean rate vanishes exponentially fast. For any strongly connected network and bounded ratios of log-likelihood functions, we obtain a lower bound on this exponential rate. For any aperiodic network we characterize the exact exponent with which probability of sample paths on which the rate of rejection deviates from the mean rate vanishes (i.e., obtain a large deviation principle) for a broader class of observation statistics including distributions with unbounded support such as Gaussian mixtures and Gamma distribution. The large deviation rate function is shown to be a function of observation model and the nodes' influences on the network as dictated by the learning rule. 
}

\noindent \textbf{Outline of the Paper.}  The rest of the paper is organized as follows. We provide the model in Section~\ref{sec:model} which defines the nodes, observations and network. This section also contains the learning rule and assumptions on model. We then provide results on rate of convergence and their proofs in Section~\ref{sec:mainresults}.  We apply our learning rule to various examples which are provided in Section~\ref{sec:performance} and some practical issues in Section~\ref{sec:prac_issues}. We conclude with a summary in Section~\ref{sec:discussion}.

\subsection{Related Work}

\add{
Literature on distributed learning, estimation and detection can divided into two broad sets. One set deals with the fusion of information observed by a group nodes at a fusion center where the communication links (between the nodes and fusion center) are either rate limited~\cite{HTcommConst1986, Han1987, DecentralizedHTGray1990, decentralizedSeqHT1993, HanAmari1994, HanAmari1998, interactiveHT2013,  optdecentralizedSeqHT2008, BinningDHT2012} or subject to channel imperfections such as fading and packet drops~\cite{channelAwareFusion2004, optLRT2005, statInfMACfading2006}. 
Our work belongs to the second set, which models the communication network as a directed graph whose vertices/nodes are agents and an edge from node $i$ to $j$ indicates that $i$ may send a message to $j$ with perfect fidelity (the link is a noiseless channel of infinite capacity). These ``protocol'' models study how message passing in a network can be used to achieve a pre-specified computational task such as distributed learning~\cite{distTrainingPreddPoor2009, learningKar2013}, general function evaluation~\cite{distSepFuncDshah}, stochastic approximations~\cite{distStocApprox2013}. Message passing protocols may be synchronous or asynchronous (such as the ``gossip'' model~\cite{randGossipAlgoDshah, ConsensusRandpathAvg2010,
ConsensusStocDis2010, ConsensusPertBlum2013,
distConsensusKar2010}).
This graphical
model of the communication, instead of a assuming a detailed physical-layer formalization, implicitly assumes a PHY/MAC-layer abstraction where sufficiently high data rates are available to send the belief vectors with desired precision when nodes are within each others'€™ communication range. A missing edge indicates the corresponding link has zero capacity. }

\add{
Due to the large body of work in distributed detection, estimation and merging of opinions, we provide a long yet detailed summary of all the related works and their relation to our setup. Readers familiar with these works can skip to Section~\ref{sec:model} without loss of continuity.
}

\add{
Several works~\cite{JadbabaieMST:12,ShahinAli:CDC2013,Jadbabaie13informationheterogeneity,KamiarRahnama:CDC2010,OlfatiSaberFES:05belief} consider an update rule which uses local Bayesian updating combined with a linear consensus strategy on the beliefs~\cite{DeGroot:74consensus} that enables all nodes in the network identify the true hypothesis. Jadbabaie et al.~\cite{JadbabaieMST:12} characterize the ``learning rate'' of the algorithm in terms of the total variational error across the network and provide an almost sure upper bound on this quantity in terms of the KL-divergences and influence vector of agents. In Corollary~\ref{cor:social_learning} we analytically show that the proposed learning rule in this paper provides a strict improvement over linear consensus strategies~\cite{JadbabaieMST:12}. Simultaneous and independent works by Shahrampour et al.~\cite{2014arXiv1409.8606S} and Nedi\'{c} et al.~\cite{2014arXiv1410.1977N} consider a similar learning rule (with a change of order in the update steps). They obtain similar convergence and concentration results under the assumption of bounded ratios of likelihood functions. Nedi\'{c} et al.~\cite{2014arXiv1410.1977N} analyze the learning rule for time-varying graphs. Theorem~\ref{thm:ldp_vector_belief} strengthens these results for static networks by providing a large deviation analysis for a broader class of likelihood functions which includes Gaussian mixtures. 
}

\add{
Rad and Tahbaz-Salehi~\cite{KamiarRahnama:CDC2010} study distributed parameter estimation using a Bayesian update rule and average consensus on the log-likelihoods similar to (\ref{eq:bayes})--(\ref{eq:estimate}). They show that the maximum of each node's belief distribution converges in probability to the true parameter under certain analytic assumptions (such as log-concavity) on the likelihood functions of the observations. Our results show almost sure convergence and concentration of the nodes' beliefs when the parameter space is discrete and the log-likelihood function is concave. Kar et al. in~\cite{distEstimationCommKar2012} consider the problem of distributed estimation of an unknown underlying parameter where the nodes make noisy observations that are non-linear functions of an unknown global parameter. They form local estimates using a quantized message-passing scheme over randomly-failing communication links, and show the local estimators are consistent and asymptotically normal. Note that for any general likelihood model, static strongly connected network and discrete parameter spaces, our Theorem~\ref{thm:expconv} strengthens the results of distributed estimation (where the error vanishes inversely with the square root of total number of observations) by showing exponentially fast convergence of the beliefs. Furthermore, Theorem~\ref{thm:expconvhoeffding} and \ref{thm:ldp_vector_belief} strengthen this by characterizing the rate of convergence. 
}
 
\add{
Similar non-Bayesian update rules have been in the context of one-shot merging of opinions~\cite{OlfatiSaberFES:05belief} and beliefs in~\cite{SaligramaDistDetect:2006} and \cite{SaligramaDistDetect:ACC2006}. Olfati-Saber et al.~\cite{OlfatiSaberFES:05belief} studied an algorithm for distributed one-shot hypothesis testing using belief propagation (BP), where nodes perform average consensus on the log-likelihoods under a single observation per node.  The nodes can achieve a consensus on the product of their local likelihoods. A benefit of our approach is that nodes do not need to know each other's likelihood functions or indeed even the space from which their observations are drawn. Saligrama et al.~\cite{SaligramaDistDetect:2006} and Alanyali et al.~\cite{SaligramaDistDetect:ACC2006}, consider a similar setup of belief propagation (after observing single event) for the problem of distributed identification of the MAP estimate (which coincides with the true hypothesis for sufficiently large number of observations) for certain balanced graphs. Each node passes messages which are composed by taking a product of the recent messages then taking a weighted average over all hypotheses. Alanyali et al.~\cite{SaligramaDistDetect:ACC2006} propose modified BP algorithms that achieves MAP consensus for arbitrary graphs. Though the structure of the message composition of the BP algorithm based message passing is similar to our proposed learning rule, we consider a dynamic setting in which observations are made infinitely often. Our rule incorporates new observation every time a node updates its belief to learn the true hypothesis. Other works study collective MAP estimation when nodes communicate discrete decisions based on Bayesian updates~\cite{HarelMST14,mueller2013general}
Harel et el. in~\cite{HarelMST14} study a two-node model where agents exchange decisions rather than beliefs and show that unidirectional transmission increases the speed of convergence over bidirectional exchange of local decisions. Mueller-Frank~\cite{mueller2013general} generalized this result to a setting in which nodes similarly exchange local strategies and local actions to make inferences. 
}

\add{
Several recently-proposed models study distributed sequential binary hypothesis testing detecting between different means with Gaussian~\cite{distDectSahu2016} and non-Gaussian observation models~\cite{distDectBajovic2012}. Jakovetic et al.~\cite{distDectBajovic2012} consider a distributed hypothesis test for i.i.d observations over time and across nodes where nodes exchange weighted sum of a local estimate from previous time instant and ratio of likelihood functions of the latest local observation with the neighbors. When the network is densely connected (for instance, a doubly stochastic weight matrix), after sufficiently long time nodes gather all the observations throughout network. By appropriately choosing a local threshold for local Neyman-Pearson test, they show that the performance of centralized Neyman-Pearson test can achieved locally. In contrast, our $M$-ary learning rule applies for observations that are correlated across nodes and exchanges more compact messages i.e., the beliefs (two finite precision real values for binary hypothesis test) as opposed to messages composed of the raw observations (in the case of $\mathbb{R}^d$ Gaussian observations with $d \gg 2$, $d$ finite precision real values for binary hypothesis test). Sahu and Kar~\cite{distDectSahu2016} consider a variant of this test for the special case of Gaussians with shifted mean and show that it minimizes the expected stopping times under each hypothesis for given detection errors.  
}

\section{The Model}
\label{sec:model}

\noindent \textbf{Notation}:  We use boldface for vectors and denote 
the $i$-th element of vector $\mathbf{v}$ by $v_{i}$. 
We let $[n]$ denote $\{1, 2, \ldots, n\}$,
$\mathcal{P}(A)$ the set of all probability distributions on a set $A$, 
$| A|$ denotes the number of elements in set $A$,
$\Ber(p)$ the Bernoulli distribution with parameter $p$, 
and $D(P_{Z}||P_Z')$ the Kullback--Leibler (KL) divergence between two probability distributions $P_Z, P_Z' \in \mc{P}(\mc{Z})$. Time is discrete and denoted by $t \in \{0, 1, 2, \ldots \}$. If $a \in A$, then $1_{a}(.) \in \mathcal{P}(A)$ denotes the probability distribution which assigns probability one to $a$ and zero probability to the rest of the elements in $A$. \add{Let $\mbf{x} \leq \mbf{y}$ denote $x_i \leq y_i$ for each $i$-th element of vector $\mbf{x}$ and $\mbf{y}$. Let $\mbf{1}$ denote the vector of where each element is $1$. For any $F \subset \mathbb{R}^{M-1}$, let $F^{o}$ be the interior of $F$ and $\bar{F}$ the closure. For $\epsilon > 0$ let $F_{\epsilon^+} = \{\mbf{x} + \delta \mbf{1}, \forall \, 0 < \delta \leq \epsilon \text{ and }\mbf{x} \in F \}$, $F_{\epsilon^-} = \{\mbf{x} - \delta \mbf{1}, \forall \, 0 < \delta \leq \epsilon \text{ and }\mbf{x} \in F\}$.
}

\subsection{Nodes and Observations}
\label{subsec:nodes_obs}
Consider a group of $n$ individual nodes. Let $\Theta = \{ \theta_1, \theta_2, \ldots, \theta_M\}$ denote a finite set of $M$ parameters which we call \textit{hypotheses}: each $\theta_i$ denotes a hypothesis. \add{At each time instant $t$, every node $i \in [n]$ makes an observation $\samp{i}{t} \in \mc{X}_i$, where $\mc{X}_i$ denotes the observation space of node $i$. The joint observation profile at any time $t$ across the network, $\{\samp{1}{t}, \samp{2}{t}, \ldots, \samp{n}{t} \}$, is denoted by $\samp{}{t} \in \mathcal{X}$, where $\mathcal{X} = \mathcal{X}_1 \times \mathcal{X}_2 \times \ldots \times \mathcal{X}_n$. The joint likelihood function for all $X \in \mathcal{X}$ given $\theta_k$ is the true hypothesis is denoted as $\dist{}{X}{\theta_k}$. We assume that the observations are statistically governed by a fixed global ``true hypothesis'' $\theta^* \in \Theta$ which is unknown to the nodes. 
Without loss of generality we assume that $\theta^* = \theta_M$. Furthermore, we assume that no node in network knows the joint likelihood functions $\{ \dist{}{\cdot}{\theta_k}\}_{k = 1}^{M}$ but every node $i \in [n]$ knows the \textit{local likelihood functions} $\{ \dist{i}{\cdot}{\theta_k}\}_{k = 1}^{M}$, where $\dist{i}{\cdot}{\theta_k}$ denotes the $i$-th marginal of $\dist{}{\cdot}{\theta_k}$.}
\remove{We consider the case where every node $i \in [n]$ is associated with a set of probability distributions $\{ \dist{i}{\cdot}{\theta} : \theta \in \Theta\}$, where $\dist{i}{\cdot}{\theta}$ is a distribution conditioned on $\theta$ being the true hypothesis.} 
Each node's observation sequence (in time) is conditionally independent and identically distributed (i.i.d) \add{but the observations might be correlated across the nodes at any given time.}

In this setting, nodes attempt to learn the ``true hypothesis'' $\theta_M$ using their knowledge of \add{$\{ \dist{i}{\cdot}{\theta_k}\}_{k = 1}^{M}$. In isolation, if $\dist{i}{\cdot}{\theta_k} \neq \dist{i}{\cdot}{\theta_M}$ for some $k \in [M-1]$, node $i$ can rule out hypothesis $\theta_k$ in favor of $\theta_M$ exponentially fast with an exponent which is equal to $\kl{\dist{i}{\cdot}{\theta_M}}{\dist{i}{\cdot}{\theta_k}}$~\cite[Section~11.7]{Cover:1991:EIT:129837}.}
\remove{$\{ \dist{i}{\cdot}{\theta} :  \theta \in \Theta \}$. It is not hard to see that if $\dist{i}{\cdot}{\theta} \neq \dist{i}{\cdot}{\theta^*}$, for some $\theta \neq \theta^*$, node $i$ can exponentially rule out hypothesis $\theta$ in favor of $\theta^*$ with an exponent which is equal to $\kl{\dist{i}{\cdot}{\theta^*}}{\dist{i}{\cdot}{\theta}}$.} Hence, for a given node the KL-divergence between the distribution of the observations conditioned over the hypotheses is a useful notion which captures the extent of distinguishability of the hypotheses. Now, define 
\begin{align*}
\bar{\Theta}_i &= \{ k \in [M] : \dist{i}{\cdot}{\theta_k} = \dist{i}{\cdot}{\theta_M}\}
	\\
	&= \{ k \in [M] : \kl{\dist{i}{\cdot}{\theta_M}}{\dist{i}{\cdot}{\theta_k}} \neq 0 \}.
\end{align*}
In other words, let $\bar{\Theta}_i$ be the set of all hypotheses that are \textit{locally indistinguishable} to node $i$.  In this work, we are interested in the case where $|\bar{\Theta}_i| > 1$ for some node $i$, but the true hypothesis $\theta_M$ is \textit{globally identifiable} (see \eqref{eq:globalID}). 

\begin{assumption}
\label{assume:distinguish}
For every pair $k \neq j$, there is at least one node $i \in [n]$ for which the KL-divergence $ \kl{ \dist{i}{\cdot}{\theta_k} }{ \dist{i}{\cdot}{\theta_j}}$ is strictly positive. 
\end{assumption}

In this case, we ask whether nodes can collectively go beyond the limitations of their local observations and learn $\theta_M$. Since
\begin{align}
\{ \theta_M \} = \bar{\Theta}_1 \cap\bar{\Theta}_2 \cap \ldots \cap \bar{\Theta}_n,
\label{eq:globalID}
\end{align}
it is straightforward to see that Assumption~\ref{assume:distinguish} is a sufficient \remove{and necessary} condition for the global identifiability of $\theta_M$ when only marginal distributions are known at the nodes. Also, note that this assumption does not require the existence of a single node that can distinguish $\theta_M$ from all other hypotheses $\theta_k$, where $k \in [M-1]$. We only require that for every pair $k \neq j$, there is at least one node $i \in [n]$ for which $\dist{i}{\cdot}{\theta_k} \neq \dist{i}{\cdot}{\theta_j}$.

\add{
Finally, we define a probability triple $\left( \Omega, \mathcal{F}, \P^{\theta_M} \right)$, where $\Omega = \{\omega: \omega = (\samp{}{0}, \samp{}{1}, \ldots), \, \forall \, \samp{}{t} \in \mathcal{X}, \, \forall \, t\}$, $\mathcal{F}$ is the $\sigma-$ algebra generated by the observations and $\P^{\theta_M}$ is the probability measure induced by paths in $\Omega$, i.e., $\P^{\theta_M} = \prod_{t = 0}^{\infty} \dist{}{\cdot}{\theta_M}$. We use $\expe^{\theta_M}[\cdot]$ to denote the expectation operator associated with measure $\P^{\theta_M}$. For simplicity we drop $\theta_M$ to denote $\P^{\theta_M}$ by $\P$
 and denote $\expe^{\theta_M}[\cdot]$ by $\expe[\cdot]$.}

\subsection{Network}
\label{subsec:network}
We model the communication network between nodes via a directed graph with vertex set $[n]$. We define the neighborhood of node $i$, denoted by $\mathcal{N}(i)$, as the set of all nodes which have an edge starting from themselves to node $i$. This means if node $j \in \mathcal{N}(i)$, it can send the information to node $i$ along this edge. In other words, the neighborhood of node $i$ denotes the set of all sources of information available to it. \add{Moreover, we assume that the nodes have knowledge of their neighbors $\mathcal{N}(i)$ only and they have no knowledge of the rest of the network~\cite{Dimakis10gossipsurvey}.}

\begin{assumption}
\label{assume:network}
The underlying graph of the network is strongly connected, i.e.\ for every $i, j \in [n]$ there exists a directed path starting from node $i$ and ending at node $j$.
\end{assumption}

We consider the case where the nodes are connected to every other node in the network by at least one multi-hop path, i.e.\ a strongly connected graph allows the information gathered to be disseminated at every node throughout the network. Hence, such a network enables learning even when some nodes in the network may not be able to distinguish the true hypothesis on their own, i.e.\ $|\bar{\Theta}_i|>1$ for some nodes.

\subsection{The Learning Rule} 
\label{subsec:learning_rule}

In this section we provide a learning rule for the nodes to learn $\theta_M$ by collaborating with each other through the local communication alone. 

We begin by defining a few variables required in order to define the learning rule. At every time instant $t$ each node $i$ maintains a \remove{estimate vector} \add{private belief vector} $\mathbf{\est{i}{t}} \in \mathcal{P}(\Theta)$ and a \add{public belief vector} $\mathbf{\bel{i}{t}} \in \mathcal{P}(\Theta)$, which are probability distributions on $\Theta$. The social interaction of the nodes is characterized by a stochastic matrix $W$. More specifically, weight $W_{ij} \in [0, 1]$ is assigned to the edge from node $j$ to node $i$ such that $W_{ij} > 0$ if and only if $j \in \mathcal{N}(i)$ and  $W_{ii} = 1 - \sum_{j = 1}^{n} W_{ij}$. The weight $W_{ij}$ denotes the confidence node $i$ has on the information it receives from node $j$. 

The steps of learning are given below. Suppose each node $i$ starts with an initial \add{private belief vector} $\mathbf{\est{i}{0}}$. At each time $t = 1, 2,\ldots$ the following events happen:
	\begin{enumerate}
	\item Each node $i$ draws a conditionally i.i.d observation $\samp{i}{t} \sim \dist{i}{\cdot}{\theta_M}$.
	\item Each node $i$ performs a local Bayesian update on $\mathbf{\est{i}{t-1}}$ to form $\mathbf{\bel{i}{t}}$ using the following rule.
	 For each $k \in [M]$, 
		\begin{align}
		\bel{i}{t}(\theta_k) = \frac{ \dist{i}{\samp{i}{t}}{\theta_k} \est{i}{t-1}(\theta_k) }{ \sum_{a \in [M]} \dist{i}{\samp{i}{t}}{\theta_a} \est{i}{t-1}(\theta_a) }.
		\label{eq:bayes}
		\end{align}
	\item Each node $i$ sends the message $\mathbf{\msg{i}{t}} = \mathbf{\bel{i}{t}}$ to all nodes $j$ for which $i \in \mathcal{N}(j)$. Similarly receives messages from its neighbors $\mathcal{N}(i)$.
	\item Each node $i$ \add{updates its private belief} of every $\theta_k$, by averaging the log beliefs it received from its neighbors. For each  $k \in [M]$,
		\begin{align}
		\est{i}{t}(\theta_k) = \frac{ \exp \left( \sum_{j = 1}^{n} W_{ij} \log \bel{j}{t}(\theta_k) \right)
			}{
			\sum_{a \in [M]} \exp \left( \sum_{j = 1}^{n} W_{ij} \log \bel{j}{t}(\theta_a) \right)
			}.
		\label{eq:estimate}
		\end{align}
	\end{enumerate}
	
Note that the \add{private belief vector} $\mathbf{\est{i}{t}}$ remain locally with the nodes while their \add{public belief vectors} $\mathbf{\bel{i}{t}}$ are exchanged with the neighbors \add{as implied by their nomenclature}.

Along with the weights, the network can be thought of as a weighted strongly connected network. Hence, from Assumption~\ref{assume:network}, we have that weight matrix $W$ is irreducible. In this context we recall the following fact.
\begin{fact}[Section~2.5 of Hoel et.\ al.~\cite{HoelPort:1972}]
\label{fact:stationarydist}
Let $W$ be the transition matrix of a Markov chain. If $W$ is irreducible then the stationary distribution of the Markov chain denoted by $\mathbf{v} = \left[v_1, v_2, \ldots, v_n\right]$ is the normalized left eigenvector of $W$ associated with eigenvalue 1 and it is given as
\begin{equation}
\label{eq:eigenvector}
v_i = \sum_{j = 1}^{n}v_j W_{ji}.
\end{equation}
Furthermore, all components of $\mathbf{v}$ are strictly positive.
If the Markov chain is aperiodic, then
\begin{equation}
\label{eq:std_dist_aperiodic}
\lim_{t \to \infty} W^{t}(i, j) = v_j, \quad i,j \in [n].
\end{equation}
If the chain is periodic with period $d$, then for each pair of states $i, j \in [n]$, there exists an integer $r \in [d]$, such that $W^{t}(i,j) = 0$ unless $t = md + r$ for some nonnegative integer $m$, and
\begin{equation}
\label{eq:std_dist_periodic}
\lim_{m \to \infty} W^{md + r}(i,j) = v_j d.
\end{equation}
\end{fact}

\noindent
In the social learning literature, the eigenvector $\mathbf{v}$ also known as the eigenvector centrality; is a measure of social influence of a node in the network. In particular we will see that $v_i$ determines the contribution of node $i$ in the collective network learning rate.

The objective of learning rule is to ensure that the \add{private belief} vector $\mathbf{\est{i}{t}}$ of each node $i \in [n]$ converges to $\mathbf{1}_{M}(\cdot)$. Note that our learning rule is such that if the initial belief of any $\theta_k, k \in [M]$, for some node is zero then beliefs of that $\theta_k$ remain zero in subsequent time intervals. Hence, we make the following assumption.

\begin{assumption}
For all $i \in [n]$, the initial \add{private belief} $\est{i}{0}(\theta_k) > 0$ for every $k \in [M]$.
\label{assume:initial_est}
\end{assumption}


\section{Main Results} 
\label{sec:mainresults}

\subsection{The Criteria for Learning}
\add{Before we present our main results, we discuss the metrics we use to evaluate the performance of a learning rule in the given distributed setup.}

\begin{definition} [\add{Rate of Rejection of Wrong Hypothesis}]
\add{For any node $i \in [n]$ and $k \in [M-1]$, define the following 
\begin{equation}
\rho_{i}^{(t)} (\theta_k) 
\defeq 
-\frac{1}{t} \log \est{i}{t}(\theta_k).
\end{equation}
The rate of rejection of $\theta_k$ in favor of $\theta_M$ at node $i$ is defined as
\begin{equation}
\rho_i (\theta_k) 
\defeq 
\liminf_{t  \to \infty} \rho_{i}^{(t)} (\theta_k).
\end{equation} 
Now, let 
\begin{align}
\estc{i} 
\defeq 
\left[\est{i}{t}(\theta_1), \est{i}{t}(\theta_2), \ldots, \est{i}{t}(\theta_{M-1})\right]^T,
\end{align}
then,
\begin{align}
\boldsymbol{\rho}_{i}^{(t)} 
\defeq
-\frac{1}{t} \log \estc{i},
\end{align}
and the rate of rejection at node $i$ is defined as 
\begin{align}
\boldsymbol{\rho}_i 
\defeq
\liminf_{t \to \infty} \boldsymbol{\rho}_{i}^{(t)} .
\end{align}}
\end{definition}
\add{If $\rho_i(\theta_k) > 0$ for all $k \in [M-1]$, under a given learning rule} the belief vectors of each node not only converge to the true hypothesis, they converge exponentially fast. Another way to measure the performance of a learning rule is the rate at which belief of true hypothesis converges to one.

\begin{definition}[\add{Rate of Convergence to True Hypothesis}]
For any $i \in [n]$ and $k \in [M-1]$, the rate of convergence to $\theta_M$, denoted by $\mu_{i}$ is defined as
\begin{equation}
\mu_{i} \defeq \add{\liminf_{t  \to \infty}}-\frac{1}{t}\log (1 - \est{i}{t}(\theta_M)).
\end{equation} 
\end{definition}

\begin{definition}[\add{Rate of Social Learning}]
The total variational error across the network when the underlying true hypothesis is $\theta_k$ \add{(where we allow the true hypothesis to vary, i.e. $\theta^* = \theta_k$ for any $k \in [M]$ instead of assuming that it is fixed at $\theta^* = \theta_M$)} is given as
\begin{equation}
e^{(t)}(k) = \frac{1}{2} \sum_{i = 1}^{n}||\est{i}{t}(\cdot) - \mathbf{1}_{k}(\cdot)|| = \sum_{i = 1}^{n} \sum_{j \neq k} \est{i}{t}(\theta_j).
\end{equation}
This equals the total probability that all nodes in the network assign to ``wrong hypotheses''. Now, define 
\begin{equation}
e^{(t)} \defeq \max_{k \in [M]} e^{(t)}(k).
\end{equation}
The rate of social learning is defined as the rate at which total variational error, $e^{(t)}$, converges to zero and mathematically it is defined as
\begin{equation}
\rho_{L} \defeq \liminf_{t \to \infty} - \frac{1}{t}\log e^{(t)}.
\end{equation}
\end{definition}
The above notion to evaluate the learning rule has been used in the social learning literature such as~\cite{Jadbabaie13informationheterogeneity}. For a given network and a given observation model for nodes, $\rho_L$ gives the least rate of learning guaranteed in the network. It is straightforward to see that with a characterization for $\rho_i (\theta_k)$ for all $k \in [M-1]$ we obtain the least rate of convergence to true hypothesis, $\mu_i$, and the least rate of social learning, $\rho_{L}$, guaranteed under a given learning rule.

\subsection{Learning: Convergence to True Hypothesis}

\begin{definition}[\add{Network Divergence}]
\add{For all $k \in [M-1]$, the network divergence between $\theta_M$ and $\theta_k$, denoted by $K(\theta_M, \theta_k)$, is defined as 
\begin{equation}
K(\theta_M, \theta_k) \defeq \sum_{i=1}^{n} v_{i} \kl{ \dist{i}{\cdot}{\theta_M} }{ \dist{i}{\cdot}{\theta_k} },
\label{eq:network_div}
\end{equation} 
$\mathbf{v} = \left[v_1, v_2, \ldots, v_n\right]$ is the normalized left eigenvector of $W$ associated with eigenvalue 1.
}
\end{definition} 
\noindent
Fact~\ref{fact:stationarydist} together with Assumption~\ref{assume:distinguish} guarantees that $K(\theta_M, \theta_k)$ is strictly positive for every $k \in [M-1]$. 

\begin{theorem}[\add{Rate of Rejecting Wrong Hypotheses, $\boldsymbol{\rho}_i$}]
\label{thm:expconv} 
Let $\theta_M$ be the true hypothesis. Under the Assumptions~\ref{assume:distinguish}--\ref{assume:initial_est}, for every node in the network, the private belief (and hence the public belief) under the proposed learning rule converges to true hypothesis exponentially fast with probability one. Furthermore, the rate of rejecting hypothesis $\theta_k$ in favor of $\theta_M$ is given by the network divergence between $\theta_M$ and $\theta_k$. Specifically, we have
\add{
\begin{align}
\lim_{t \to \infty} \mathbf{q_i^{(t)}} = \mbf{1}_M \quad \P\text{-a.s.}
\end{align}
and 
\begin{align}
\boldsymbol{\rho}_i = - \lim_{t \to \infty} \frac{1}{t} \log \estc{i} =  \mbf{K} \quad \P\text{-a.s.}
\end{align}
where
\begin{align}
\mbf{K} = \left[K(\theta_M, \theta_1), K(\theta_M, \theta_2), \ldots, K(\theta_M, \theta_{M-1}) \right]^T.
\end{align}}
\end{theorem}


Theorem~\ref{thm:expconv} establishes that the \add{beliefs} \remove{estimates}  of wrong hypotheses, $\theta_k$ for $k \in [M-1]$, vanish exponentially fast and it characterizes the exponent with which a node rejects $\theta_k$ in favor of $\theta_M$. This rate of rejection is a function of the node's ability to distinguish between the hypotheses given by the KL-divergences and structure of the weighted network which is captured by the eigenvector centrality of the nodes. Hence, every node influences the rate in two ways. Firstly, if the node has higher eigenvector centrality (\emph{i.e.} the node is centrality located), it has larger influence over the \add{beliefs} \remove{estimates} of other nodes as a result has a greater influence over the rate of exponential decay as well. Secondly, if the node has high KL-divergence (\emph{i.e} highly informative observations that can distinguish between $\theta_k$ and $\theta_M$), then again it increases the rate. If an influential node has highly informative observations then it boosts the rate of rejecting $\theta_k$ by improving the rate. We will illustrate this through a few numerical examples in Section~\ref{subsec:rate_influence}.

We obtain lower bound on the rate of convergence to the true hypothesis and rate of learning as corollaries to Theorem~\ref{thm:expconv}.

\begin{corollary}[Lower Bound on Rate of Convergence to $\theta_M$]
\label{cor:conv}
Let $\theta_M$ be the true hypothesis. Under the Assumptions~\ref{assume:distinguish}--\ref{assume:initial_est}, for every $i \in [n]$, the rate of convergence to $\theta_M$ can be lowered bounded as 
\begin{equation}
\mu_i \geq   \min_{k \in [M-1]} K(\theta_M, \theta_k)  \quad \P\mbox{-a.s.}
\end{equation}
\end{corollary}

\begin{corollary} [Lower Bound on Rate of Learning]
\label{cor:social_learning}
Let $\theta_M$ be the true hypothesis. Under the Assumptions~\ref{assume:distinguish}--\ref{assume:initial_est}, the rate of learning $\rho_L$ across the network is lower bounded by,
\begin{align*}
\rho_L &\geq   \min_{i,j \in [M]}{ K(\theta_i, \theta_j)}  \quad \P\mbox{-a.s.}
\end{align*}
\end{corollary}

\begin{remarks}
\label{rem:lambda}
Jadbabaie et.\ al. proposed a learning rule in~\cite{JadbabaieMST:12}, which differs from the proposed rule at the private belief vector $\mathbf{\est{i}{t}}$ formation step. Instead of averaging the log beliefs, nodes average the beliefs received as messages from their neighbors. In~\cite{Jadbabaie13informationheterogeneity}, Jadbabaie et.\ al. provide an upper bound on the rate of learning $\rho_L$ obtained using their algorithm. They show
\begin{equation}
\rho_L \leq \alpha \min_{i,j \in [M]}{ K(\theta_i, \theta_j)}  \quad \P\mbox{-a.s.}
\label{eq:lambda_upp_bound}
\end{equation}
where $\alpha$ is a constant strictly less than one. Corollary~\ref{cor:social_learning} shows that lower bound on $\rho_L$ using the proposed algorithm is greater than the upper bound provided in equation~\ref{eq:lambda_upp_bound}.
\end{remarks}

\remove{Shahrampour and Jadbabaie~[2] for a doubly stochastic $W$ provide a closed-form of an lower bound on $\mu_i$, where nodes communicate according to a random gossip scheme followed by a stochastic optimization based (high complexity) update performed by all nodes. Lower bound is given by}
\remove{We show that using the proposed learning rule which has low complexity update steps, we achieve the same lower bound on $\mu_i$ as obtained in~[2] whenever the weight matrix $W$ is doubly stochastic. In a recent independent study~[11] et al. propose an extension of their distributed stochastic optimization based learning rule in~[11] and its closed form solution is similar to the proposed learning rule. Both works are applicable to aperiodic networks and hence they achieve the similar convergence results for aperiodic networks.}

\subsection{Concentration under Bounded Log-likelihood ratios}

\add{Under very mild assumptions, Theorem~\ref{thm:expconv} shows that the belief of a wrong hypothesis $\theta_k$ for $k \in [M-1]$ converging to zero exponentially fast at rate equal to the network divergence between $\theta_M$ and $\theta_k$, $ K(\theta_M, \theta_k)$, with probability one. We strength this result under the following assumption.}

\begin{assumption}
\label{assume:boundedMG}
There exists a positive constant $L$ such that
	\begin{equation}
	\max_{i \in [n]}\max_{j,k \in [M]}\sup_{X \in \mc{X}_i} \left| \log \frac{ \dist{i}{X}{\theta_j}}{ \dist{i}{X}{\theta_k}} \right|\leq L.
	\end{equation}
\end{assumption}

\begin{theorem}[Concentration of Rate of Rejecting Wrong Hypotheses, $\rho_i^{(t)}(\theta_k)$]
\label{thm:expconvhoeffding} 
Let $\theta_M$ be the true hypothesis. Under Assumptions~\ref{assume:distinguish}--\ref{assume:boundedMG}, for every node $i \in [n]$, $k \in [M-1]$, and for all $\epsilon > 0$ we have
\begin{equation}
\lim_{t \to \infty} \frac{1}{t} \log \P \left(\rho_{i}^{(t)}(\theta_k) \leq   K(\theta_M, \theta_k) - \epsilon \right) \leq -\frac{\epsilon^2}{2 L^2d}.
 \end{equation} 
For $0< \epsilon \leq   L -   K(\theta_M, \theta_k)$, we have
\begin{align}
\lim_{t \to \infty}\frac{1}{t} \log \P \left( \rho_{i}^{(t)}(\theta_k) \geq   K(\theta_M, \theta_k) + \epsilon \right)
\nonumber
\\
\leq -\frac{1}{2 L^2d} \min \left\{\epsilon^2,   \min_{j \in [M-1]} K^2(\theta_M, \theta_j) \right\}.
\end{align} 
\add{
For $\epsilon \geq   L -   K(\theta_M, \theta_k)$ we have
\begin{align}
\lim_{t \to \infty}\frac{1}{t} \log 
\P \left( \rho_{i}^{(t)}(\theta_k) \geq   K(\theta_M, \theta_k) + \epsilon \right)
\nonumber
\\
\leq -\min_{k \in [M-1]} \left\{\frac{K(\theta_M, \theta_k)^2}{2 L^2d}\right\}.
\end{align} 
}
\end{theorem}
 
\begin{corollary}[Rate of convergence to True Hypothesis]
\label{cor:concentration_learning_rate}
Let $\theta_M$ be the true hypothesis. Under Assumptions~\ref{assume:distinguish}--\ref{assume:boundedMG}, for every $i \in [n]$, we have
\begin{equation*}
\mu_i =   \min_{k \in [M-1]} K(\theta_M, \theta_k)  \quad \P\mbox{-a.s.}
\end{equation*}
\end{corollary}

From Theorem~\ref{thm:expconv} we know that $\rho_{i}^{(t)}(\theta_k)$ converges to $  K(\theta_M, \theta_k)$ almost surely. Theorem~\ref{thm:expconvhoeffding} strengthens Theorem~\ref{thm:expconv} by showing that the probability of sample paths where $\rho_{i}^{(t)}(\theta_k)$ deviates by some fixed $\epsilon$ from $  K(\theta_M, \theta_k)$, vanishes exponentially fast. This implies that $\rho_{i}^{(t)}(\theta_k)$ converges to $  K(\theta_M, \theta_k)$ exponentially fast in probability. Also, Theorem~\ref{thm:expconvhoeffding} characterizes a lower bound on the exponent with the probability of such events vanishes and shows that periodicity of the network reduces the exponent.

\subsection{Large Deviation Analysis}

\add{
\begin{assumption}
\label{assume:finite_lmgf}
For every pair $\theta_i \neq \theta_j$ and every node $k \in [n]$, the random variable $ \left| \log \frac{\dist{k}{X_k}{\theta_i}}{\dist{k}{X_k}{\theta_j}} \right|$ has finite log moment generating function under distribution $\dist{k}{\cdot}{\theta_j}$.   
\end{assumption}
}

\add{
This is a technical assumption that it relaxes the assumption of bounded ratios of the likelihood functions in prior work~\cite{DHT_ISIT2014,LalithaJ:15allerton,2014arXiv1409.8606S,Nedic2015arXiv}. Next, we provide families of distributions which satisfy Assumption~\ref{assume:finite_lmgf} but violate Assumption~\ref{assume:boundedMG}.\\
}

\add{
\begin{remarks}
Distributions $f(X;\theta_i)$ and $f(X;\theta_j)$ for $i \neq j$ with the following properties for some positive constants $C$ and $\beta$,
\begin{align}
\label{eq:suff_assp_2}
\P_i\left(\frac{f(X;\theta_j)}{f(X;\theta_i)} \geq x \right) \leq \frac{C}{x^\beta},
\quad
\P_i\left(\frac{f(X;\theta_i)}{f(X;\theta_j)} \geq x \right) \leq \frac{C}{x^\beta},
\end{align}
satisfy Assumption~\ref{assume:finite_lmgf}. Note that~(\ref{eq:suff_assp_2}) is a sufficient condition but not a necessary condition. Examples~\ref{ex:gaussian}--\ref{ex:gamma} below do not satisfy~(\ref{eq:suff_assp_2}) yet satisfy Assumption~\ref{assume:finite_lmgf}.
\end{remarks}
}

\add{
\begin{example}[Gaussian distribution and Mixtures]
\label{ex:gaussian}
Let $f(X; \theta_1) = \mathcal{N}(\mu_1, \sigma)$ and $f(X; \theta_2) = \mathcal{N}(\mu_2, \sigma)$, then 
\begin{align}
\left| \log \frac{f(x; \theta_1)}{f(x; \theta_2)} \right|
\leq
c_1 |x| + c_2,
\end{align}
where $c_1 = \left| \frac{\mu_1 - \mu_2}{\sigma^2}\right|$ and $c_2 = \left| \frac{\mu_1^2 - \mu_2^2}{2\sigma^2}\right|$. Hence, for $\lambda \geq 0$ we have
\begin{align}
\expe \left[ e^{\lambda \left| \log \frac{f(X; \theta_1)}{f(X; \theta_2)} \right|}\right]
\leq 
e^{c_2 \lambda} \expe \left[ e^{c_1 \lambda |x|}\right]
< \infty.
\end{align}
More generally for $i\in \{1,2\}$, and $p \in [0,1]$, let
\begin{align}
f(x;\theta_i)
&= \frac{p}{\sigma \sqrt{2 \pi}}\exp\left(\frac{-(x-\alpha_i)^2}{2 \sigma^2}\right)
\nonumber
\\
&+
\frac{1-p}{\sigma \sqrt{2 \pi}} \exp\left(\frac{-(x-\beta_i)^2}{2 \sigma^2}\right).
\end{align}
Then the log moment generating function of $\left| \log \frac{f(X;\theta_1)}{f(X;\theta_2)}\right|$ is finite for all $\lambda \geq 0$.
\end{example}
}

\add{
\begin{example}[Gamma distribution]
\label{ex:gamma}
Let $f(X; \theta_1) = \frac{\beta^{\alpha_1}}{\Gamma(\alpha_1)} x^{\alpha_1 - 1} e^{-\beta x}$ and $f(X; \theta_2) = \frac{\beta^{\alpha_2}}{\Gamma(\alpha_2)} x^{\alpha_2 - 1} e^{-\beta x}$, then 
\begin{align}
\left| \log \frac{f(x; \theta_1)}{f(x; \theta_2)} \right|
\leq
c_1 |\log x| + c_2,
\end{align}
where $c_1 = \left| \alpha_1 - \alpha_2 \right|$ and $c_2 = \left| (\alpha_1 - \alpha_2)\log \beta + \log \frac{\Gamma(\alpha_2)}{\Gamma(\alpha_1)}\right|$. Hence, for $\lambda \geq 0$ we have
\begin{align}
\expe \left[ e^{\lambda \left| \log \frac{f(X; \theta_1)}{f(X; \theta_2)} \right|}\right]
\leq 
e^{c_2 \lambda} \expe \left[ e^{c_1 \lambda |\log x|}\right]
< \infty.
\end{align}
\end{example}
}

\add{
The above examples show that Assumption~\ref{assume:finite_lmgf} is satisfied for distributions which have unbounded support. In order to analyze the concentration of $\boldsymbol{\rho}_i^{(t)}$ under Assumption~\ref{assume:finite_lmgf} we replace Assumption~\ref{assume:network} with the following assumption.

\noindent
\textbf{Assumption $\mathbf{2^{\prime}}$.} The underlying graph of the network is strongly connected and aperiodic.
}

\add{
Now we provide few more definitions. Let 
\begin{align}
\mathbf{Y}^{(t)}
\defeq 
  \sum_{k = 1}^{M-1} \langle \mbf{v}, \mathbf{L}^{(t)}(\theta_k) \rangle,
\label{eq:llr_corr}
\end{align}
where $\mathbf{L}^{(t)}(\theta_k)$ is the vector of log likelihood ratios given by
\begin{align}
&\mathbf{L}^{(t)}(\theta_k)
\nonumber
\\
&= 
\left[
\log \frac{\dist{1}{\samp{1}{t}}{\theta_k}}{\dist{1}{\samp{1}{t}}{\theta_M}}, \ldots,
\log \frac{\dist{n}{\samp{n}{t}}{\theta_{k}}}{\dist{n}{\samp{n}{t}}{\theta_M}}
\right]^T.
\label{eq:llr_vector}
\end{align}
}

\add{
\begin{definition} [Moment Generating Function]
For every $\lambda_k \in \mathbb{R}$, let $\Lambda_k(\lambda_k)$ denote the log moment generating function of $\langle \mbf{v},   \mathbf{L}(\theta_k) \rangle$ given by
\begin{align}
\Lambda_k(\lambda_k) 
&\defeq
\log \expe [e^{  \lambda_k \langle \mbf{v}, \mbf{L}(\theta_k)\rangle}]
\nonumber
\\
&=
\sum_{j = 1}^{n} 
\log \expe 
\left[ 
\left\{ 
\frac{\dist{j}{X_j}{\theta_k}}{\dist{j}{X_j}{\theta_M}}  \right\}^{  \lambda_k  v_j}
\right].
\end{align}
For every $\boldsymbol{\lambda} \in \mathbb{R}^{M-1}$, let $\Lambda(\boldsymbol{\lambda})$ denote the log moment generating function of $\mathbf{Y}$ given by
\begin{align}
\Lambda(\boldsymbol{\lambda}) 
&\defeq
\log \expe [e^{\langle \boldsymbol{\lambda}, \mbf{Y}\rangle}]
=
\sum_{k = 1}^{M-1} \Lambda_k(\lambda_k).
\end{align}
\end{definition}
}

\add{
\begin{definition}[Large Deviation Rate Function]
For all $x \in \mathbb{R}$, let $I_k(x)$ denote the Fenchel-Legendre transform of $\Lambda_k(\cdot)$ and is given by
\begin{align}
\label{eq:rate_func_single_theta}
I_k(x) 
\defeq
\sup_{\lambda_k \in \mathbb{R}} \left\{ \lambda x - \Lambda_k(\lambda_k)\right\}.
\end{align}
For all $\mbf{x} \in \mathbb{R}^{M-1}$, let $I(\mbf{x})$ denote the Fenchel-Legendre transform of $\Lambda(\cdot)$ and is given by
\begin{align}
\label{eq:rate_func}
I(\mathbf{x}) 
\defeq
\sup_{\boldsymbol{\lambda} \in \mathbb{R}^{M-1}} \left\{ \left\langle \boldsymbol{\lambda}, 
\mathbf{x}
\right\rangle - \Lambda(\boldsymbol{\lambda})\right\}.
\end{align}
\end{definition}
}

\add{
\begin{theorem}[Large Deviations of $\boldsymbol{\rho}_{i}^{(t)}$]
\label{thm:ldp_vector_belief}
Let $\theta_M$ be the true hypothesis. Under Assumptions~\ref{assume:distinguish}, $2^{\prime}$, \ref{assume:initial_est}, \ref{assume:finite_lmgf}, the rate of rejection $\boldsymbol{\rho}_{i}^{(t)}$ satisfies an Large Deviation Principle with rate function $J(\cdot)$, i.e., for any set $F \subset \mathbb{R}^{M-1}$ we have
\begin{align}
\liminf_{t \to \infty} \frac{1}{t} \log 
\P \left(
 \boldsymbol{\rho}_{i}^{(t)} \in F
\right)
\geq - \inf_{\mbf{y} \in F^o} J(\mbf{y}),
\end{align}
and
\begin{align}
\limsup_{t \to \infty} \frac{1}{t} \log 
\P \left(
 \boldsymbol{\rho}_{i}^{(t)} \in F
\right)
\leq - \inf_{\mbf{y} \in \bar{F}} J(\mbf{y}),
\end{align}
where large deviation rate function $J(\cdot)$ is defined as
\begin{align}
J(\mbf{y}) 
\defeq \inf_{\mathbf{x} \in \mathbb{R}^{M-1}: g(\mathbf{x}) =\mbf{y}} 
I(\mbf{x}),
\,\, \forall \, \mbf{y} \in \mathbb{R}^{M-1},
\end{align}
where $g:\mathbb{R}^{M-1} \to \mathbb{R}^{M-1}$ is a continuous mapping given by
\begin{align}
g(\mathbf{x}) 
\defeq 
\left[
g_1(\mbf{x}), g_2(\mbf{x}), \ldots, g_{M-1}(\mbf{x})
\right]^T, 
\end{align}
and 
\begin{align}
g_k(\mathbf{x}) &\defeq x_k - \max\{0, x_1, x_2, \ldots, x_{M-1}\}.
\end{align}
\end{theorem}
}

\add{
Theorem~\ref{thm:ldp_vector_belief} characterizes the asymptotic rate of concentration
of $\boldsymbol{\rho}_{i}^{(t)}$ in any set $F \subset \mathbb{R}^{M-1}$. In other words, it characterizes the rate at which the probability of deviations in each $\rho_i^{(t)}(\theta_k)$ from the rate of rejection $K(\theta_M, \theta_k)$ for every $\theta_k$ for every $k \in [M-1]$ vanish simultaneously. It characterizes the asymptotic rate as a function of the observation model of each node (not just the bound $L$ on the ratios of log-likelihood function) and as a function of eigenvector centrality $\mbf v$. The following corollary specializes this result to obtain the individual rate of rejecting a wrong
hypothesis at every node.
}

\add{
\begin{corollary}
\label{coll:ldp_single_belief}
Let $\theta_M$ be the true hypothesis. Under Assumptions~\ref{assume:distinguish}, $2^{\prime}$, \ref{assume:initial_est}, \ref{assume:finite_lmgf}, for $0 < \epsilon \leq   K(\theta_M, \theta_k)$, $k \in [M-1]$, we have
\begin{align}
\label{eq:LDPabovemean}
\lim_{t \to \infty} \frac{1}{t} \log 
\P \left( 
\rho_i^{(t)}(\theta_k) \leq   K(\theta_M, \theta_k) - \epsilon \right)
\nonumber
\\
= 
- I_k \left(  K(\theta_M, \theta_k) - \epsilon \right),
\end{align}
and for $\epsilon > 0$, we have
\begin{align}
\label{eq:LDPbelowmean}
\lim_{t \to \infty} \frac{1}{t} \log 
\P \left( 
\rho_i^{(t)}(\theta_k) \geq   K(\theta_M, \theta_k) + \epsilon \right)
\nonumber
\\
= 
- I_k \left(   K(\theta_M, \theta_k) + \epsilon \right).
\end{align}
\end{corollary}
}

\add{
Using Theorem~\ref{thm:ldp_vector_belief} and Hoeffding's Lemma, we obtain the following corollary.
\begin{corollary} 
\label{coll:LDPtoHoeffding}
Suppose Assumption~\ref{assume:boundedMG} is satisfied for some finite $L \in \mathbb{R}$. Then for small $\epsilon $ as specified in Theorem~\ref{thm:ldp_vector_belief}, we recover the exponents of Theorem~\ref{thm:expconvhoeffding} under aperiodic networks, given by
\begin{align}
\lim_{t \to \infty} \frac{1}{t} \log \P \left( \rho_i^{(t)}(\theta_k) \geq   K(\theta_M, \theta_k) + \epsilon \right) 
\leq -\frac{\epsilon^2}{2   L^2},
\end{align}
and
\begin{align}
\lim_{t \to \infty} \frac{1}{t} \log \P \left( \rho_i^{(t)}(\theta_k) \leq   K(\theta_M, \theta_k) - \epsilon \right) 
\leq -\frac{\epsilon^2}{2   L^2}.
\end{align}
\end{corollary}
}

\add{
\begin{remarks}
Under Assumption~\ref{assume:boundedMG}, Corollary~\ref{coll:LDPtoHoeffding} shows that lower bound on the asymptotic rate of concentration
of $\boldsymbol{\rho}_{i}^{(t)}$ as characterized by Theorem~\ref{thm:expconvhoeffding} is loose in comparision to that obtained from Theorem~\ref{thm:ldp_vector_belief}.  Nedic et al.\cite{2014arXiv1410.1977N} and Shahrampour et al. \cite{2014arXiv1409.8606S} provide non-asymptotic lower bounds on the rate of concentration
of $\boldsymbol{\rho}_{i}^{(t)}$ whose asymptotic form coincides with the lower bound on rate characterized by Theorem~\ref{thm:expconvhoeffding} for  aperiodic networks. This implies that under Assumption~\ref{assume:boundedMG} Theorem~\ref{thm:ldp_vector_belief} provides a tighter asymptotic rate than  that in\cite{2014arXiv1410.1977N} and \cite{2014arXiv1409.8606S}. Hence, Theorem~\ref{thm:ldp_vector_belief} strengthens Theorem~\ref{thm:expconvhoeffding} by extending the large deviation to larger class of distributions and by capturing the complete effect of nodes' influence in the network and the local observation statistics.
\end{remarks}
}

\section{Examples}\label{sec:performance}

In this section through numerical examples we illustrate how nodes learn using the proposed scheme and examine the factors which affect the rate of rejection of wrong hypotheses and its rate of concentration.

\subsection{Factors influencing Convergence}
\label{subsec:rate_influence}

\begin{example}
\label{ex:2node}
\add{Consider a group of two nodes as shown in Figure \ref{fig:grid}, where the set of hypotheses is $\Theta = \{\theta_1, \theta_2, \theta_3, \theta_4\}$ and true hypothesis $\theta^* = \theta_4$. Observations at each node at time $t$, $X_i^{(t)}$, take values in $\mathbb{R}^{100}$ and have a Gaussian distribution. For node 1, $\dist{1}{\cdot}{\theta_1} = \dist{1}{\cdot}{\theta_3}  =  \mathcal{N}(\boldsymbol{\mu}_{11}, \boldsymbol{\Sigma})$ and $\dist{1}{\cdot}{\theta_2} = \dist{1}{\cdot}{\theta_4} =  \mathcal{N}(\boldsymbol{\mu}_{12}, \boldsymbol{\Sigma})$, and for node 2,  $\dist{2}{\cdot}{\theta_1} = \dist{2}{\cdot}{\theta_2} =  \mathcal{N}(\boldsymbol{\mu}_{21}, \boldsymbol{\Sigma})$ and $\dist{2}{\cdot}{\theta_3} = \dist{2}{\cdot}{\theta_4}  =  \mathcal{N}(\boldsymbol{\mu}_{22}, \boldsymbol{\Sigma})$, where $\boldsymbol{\mu}_{11}, \boldsymbol{\mu}_{12}, \boldsymbol{\mu}_{21}, \boldsymbol{\mu}_{22} \in \mathbb{R}^{100}$ and $\boldsymbol{\Sigma}$ is a positive semi-definite matrix of size 100-by-100. Here, node 1 can identify the column containing $\theta_4$, and node 2 can identify the row. In other words, $\bar{\Theta}_1 = \{ \theta_2, \theta_4\}$ and $\bar{\Theta}_2 = \{ \theta_3, \theta_4\}$. Also, $\theta_4 = \bar{\Theta}_1 \cap \bar{\Theta}_2$, hence $\theta_4$ is globally identifiable. }
\end{example}

\subsubsection{Strong Connectivity}
\begin{figure} [h!] 
  \centering
    \includegraphics[width=0.45\textwidth]{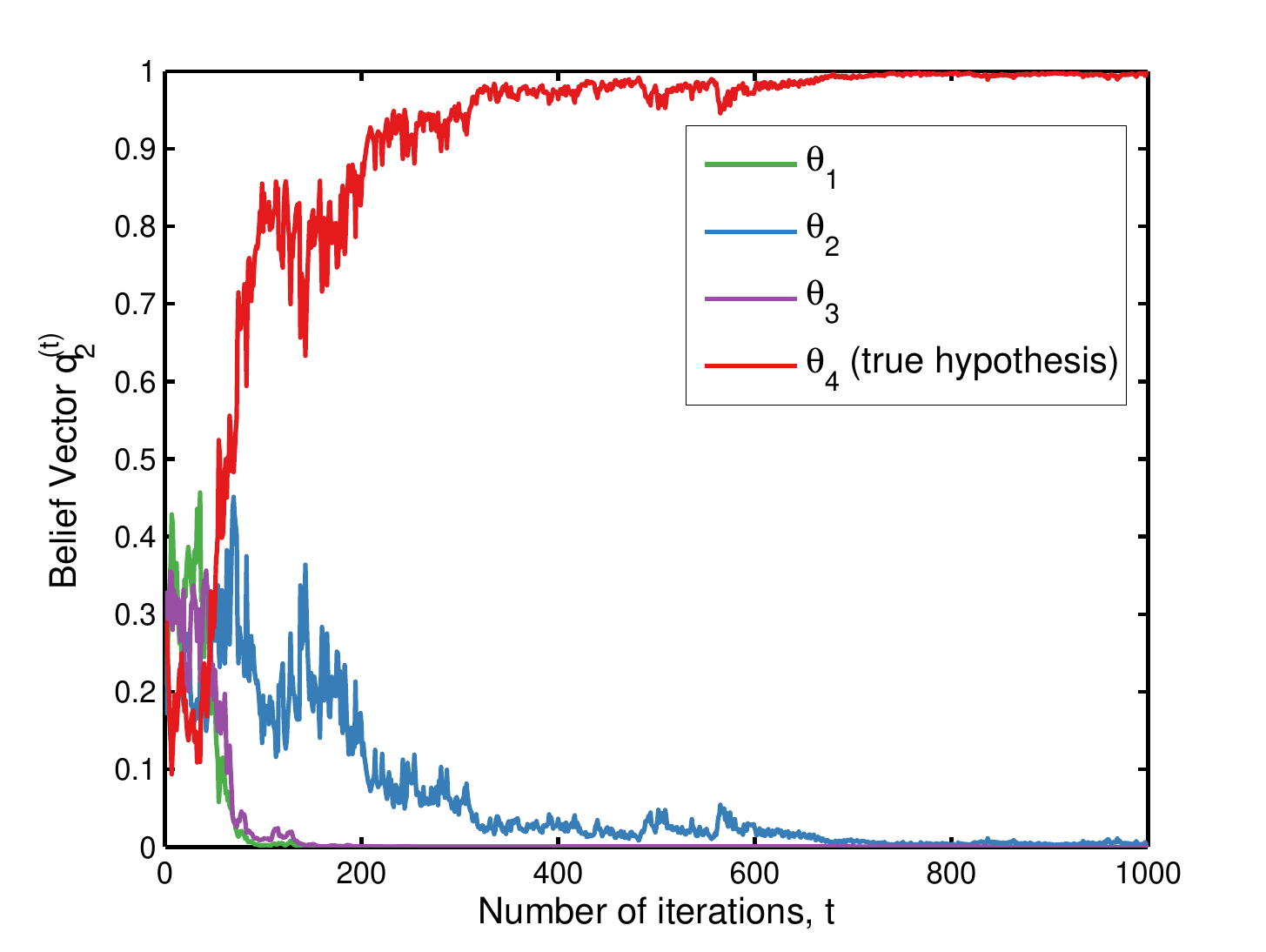}
    \caption{For the set of nodes described in Figure~\ref{fig:grid}, this figure shows the evolution of beliefs for one instance using the proposed learning rule. Belief of the true hypothesis $\theta_4$ of node 2 converges to 1 and beliefs of all other hypotheses go to zero.}
    \label{fig:two_node_4para}
\end{figure}

\begin{figure} [h!] 
  \centering
    \includegraphics[width=0.45\textwidth]{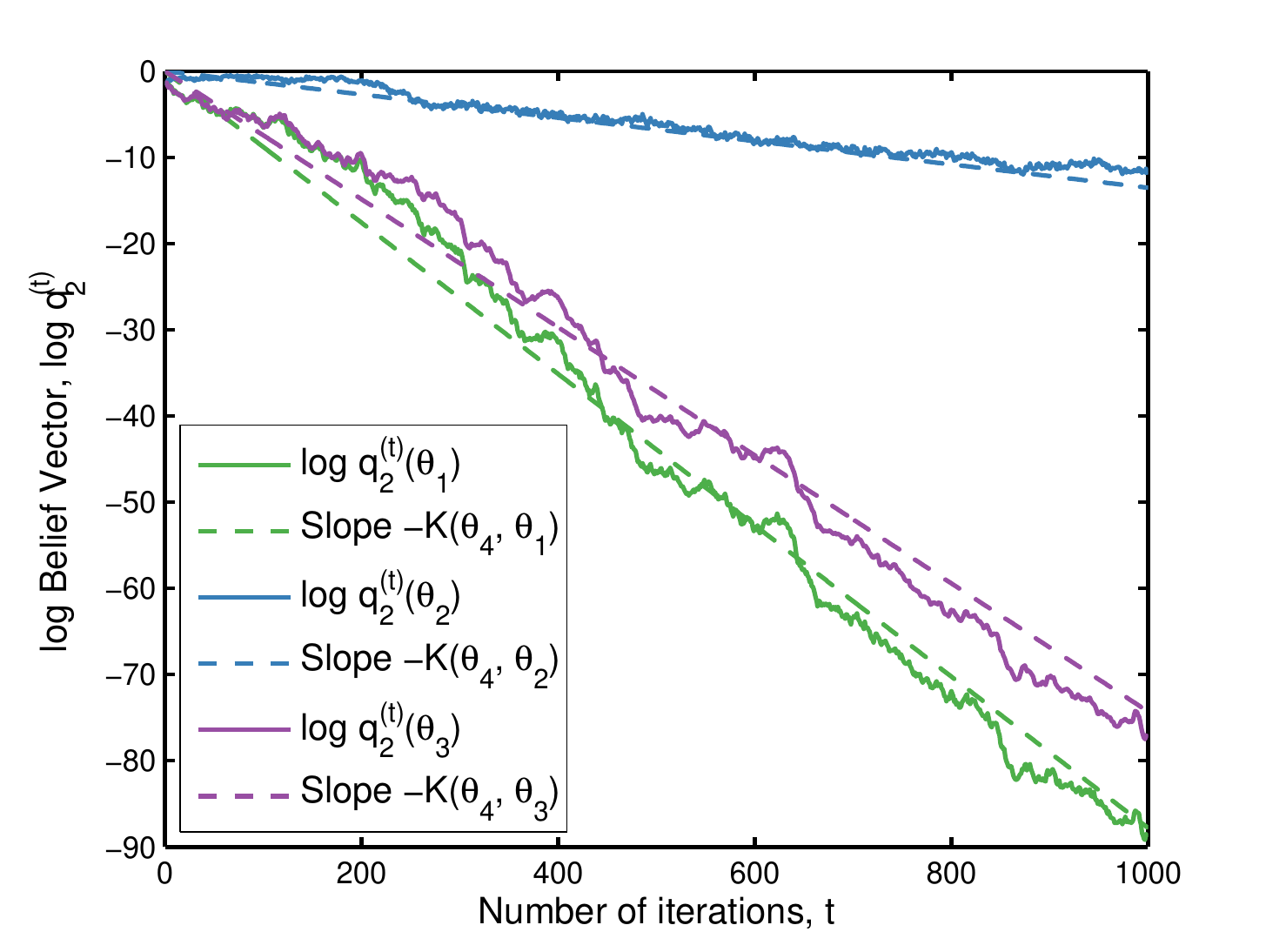}
    \caption{Figure shows the exponential decay of beliefs of $\theta_1$, $\theta_2$ and $\theta_3$ of node 2 using the learning rule.}
    \label{fig:rej_rates}
\end{figure}

Nodes are connected to each other in a network and the weight matrix is given by
\begin{align}
\label{eq:aperiodic_w}
	W = \left( \begin{array}{cc}
	0.9 & 0.1\\
	0.4 & 0.6\end{array}\right).
\end{align}
Figure~\ref{fig:two_node_4para} shows the evolution of beliefs with time for node 2 on a single sample path. We see that using the proposed learning rule, belief of $\theta_4$  goes to one while the beliefs of other hypotheses go to zero. This example shows that each node by collaboration is able to see new information which was not available through its local observations alone and both nodes learn $\theta_4$. Figure~\ref{fig:rej_rates} shows the rate of rejection of  wrong hypotheses. We see that the rate of rejection $\theta_k$ for $k \in \{1,2, 3\}$ closely follows the asymptotic rate $K(\theta_4, \theta_k)$.

\begin{figure} [h!] 
  \centering
    \includegraphics[width=0.45\textwidth]{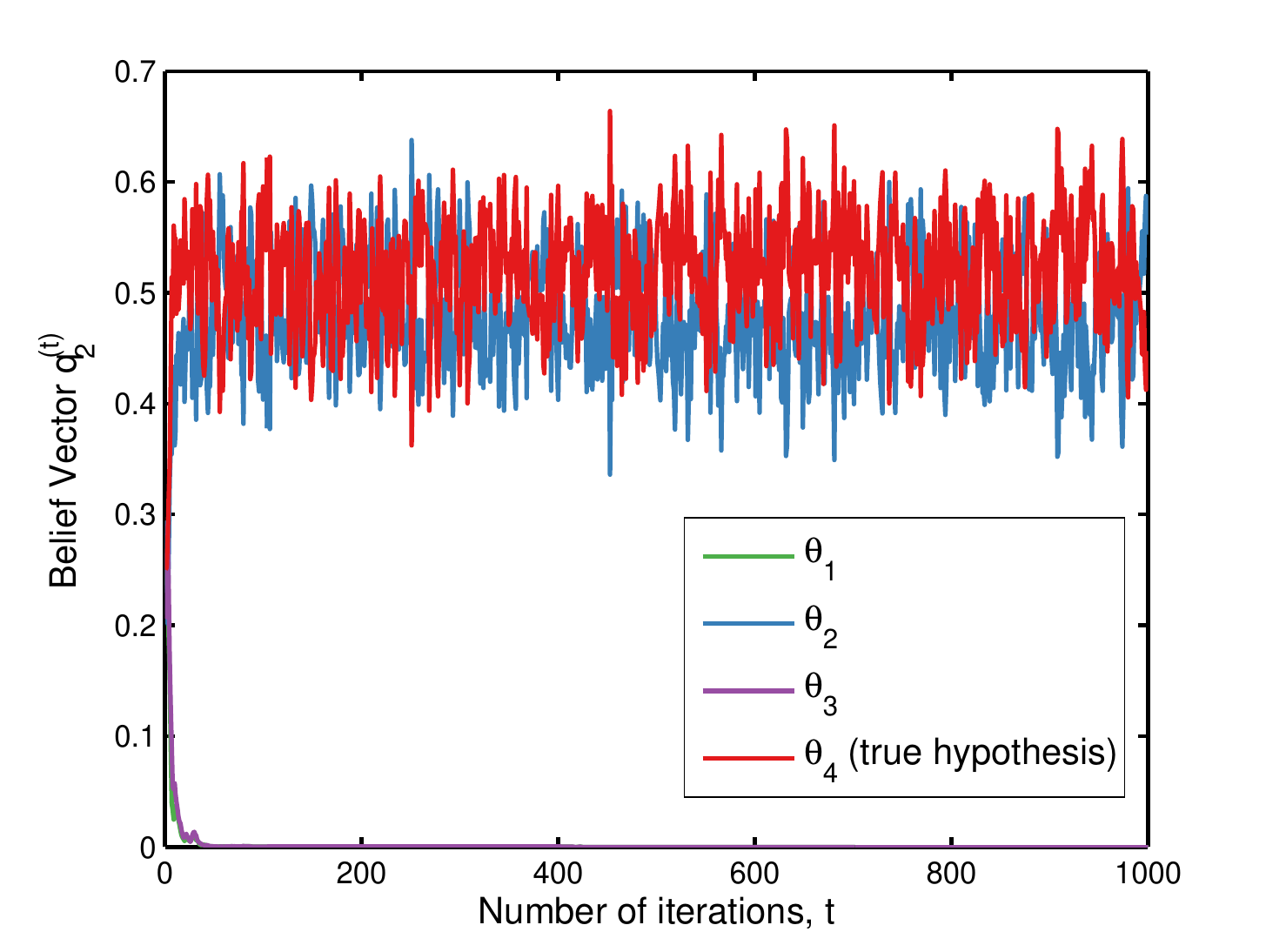}
    \caption{Figure shows the beliefs of node 2 shown in Figure~\ref{fig:grid}. When the network is not strongly connected node 2 cannot learn $\theta_4$.}
    \label{fig:2node_noNW}
\end{figure}

Suppose the nodes are connected to each other in a network which is not strongly connected and its weight matrix is given by  
\begin{align}
\label{eq:not_conn_w}
	W = \left( \begin{array}{cc}
	1 & 0\\
	0.5 & 0.5\end{array}\right).
\end{align}
Since there is no path from node 2 to node 1, the network is not strongly connected anymore. Node 2 as seen in Figure~\ref{fig:2node_noNW} does not converge to $\theta_4$. Even though node 1 cannot distinguish the elements of $\bar{\Theta}_1$ from $\theta_4$, it rejects the hypotheses in $\{\theta_1, \theta_3 \}$ in favor of $\theta_4$. This forces node 2 also to reject the set $\{\theta_1, \theta_3\}$. For node 1, $\theta_2$ and $\theta_4$ are observationally equivalent, hence their respective beliefs equal half. But node 2 oscillates between $\theta_2$ and $\theta_4$ and is unable to learn $\theta_4$. Hence, when the network is not strongly connected both nodes fail to learn.

\begin{figure} [h!] 
  \centering
    \includegraphics[width=0.45\textwidth]{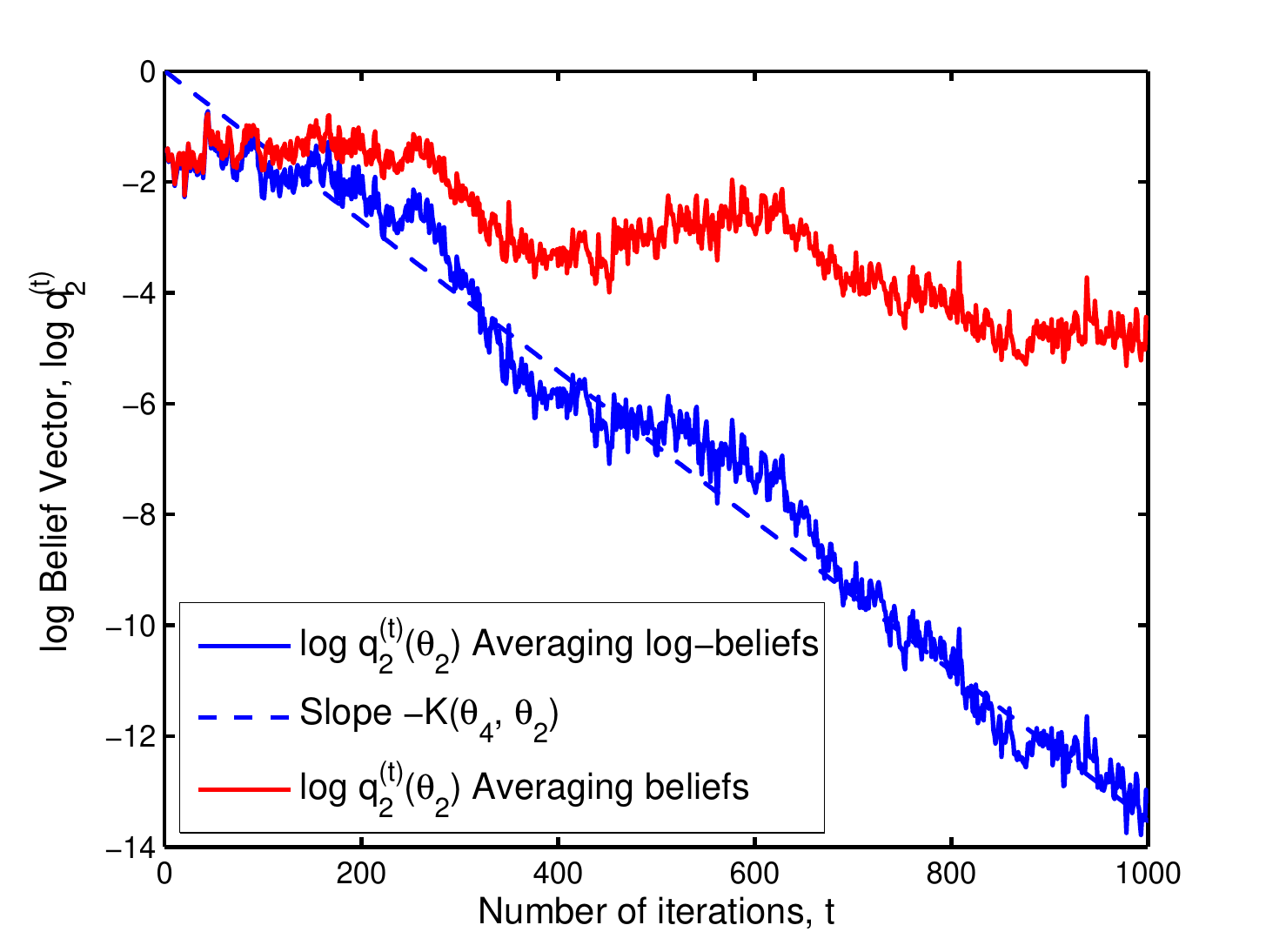}
    \caption{Figure shows that the rate of rejection of $\theta_2$ using the proposed learning rule (averaging the log beliefs) is greater than the rate of rejection of $\theta_2$ obtained using the learning rule in~\cite{JadbabaieMST:12} (averaging the beliefs).}
    \label{fig:comp_rules}
\end{figure}

In this setup we apply the learning rule considered in~\cite{JadbabaieMST:12}, where in the consensus step public beliefs are updated by averaging the beliefs received from the neighbors instead of averaging the logarithm of the beliefs. As seen in Figure \ref{fig:comp_rules}, rate of rejecting learning using the proposed learning rule is greater than the upper bound on learning rule in~\cite{JadbabaieMST:12}. Note that the precision of the belief vectors in the simulations is 8 bytes i.e. 64 bits per hypothesis. This implies the nodes each send 32 bytes per unit time, which is less than the case when nodes exchange raw Gaussian observations which may require data rate as high as 800 bytes per observation when each dimension of the Gaussian is independent.

\subsubsection{Periodicity} 
\begin{figure} [h!] 
  \centering
    \includegraphics[width=0.45\textwidth]{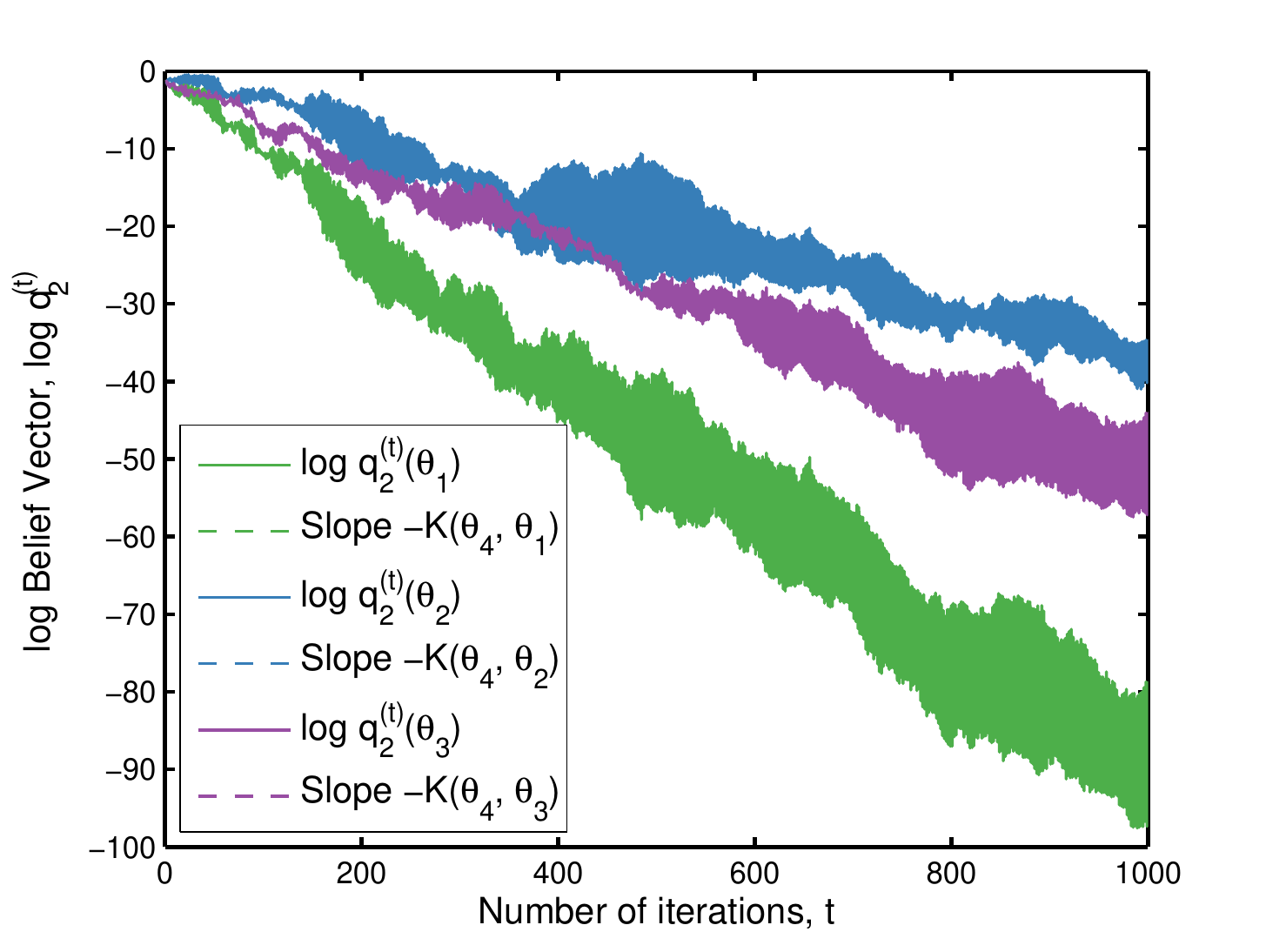}
    \caption{Figure shows the exponential decay of beliefs of $\theta_1$, $\theta_2$, and $\theta_3$ of node 2 connected to node 1 in a periodic network with period 2.}
    \label{fig:rej_rates_periodic}
\end{figure}

Now suppose the nodes are connected to each other in periodic network with period 2 and the weight matrix given by
\begin{align}
\label{eq:periodic_w}
	W = \left( \begin{array}{cc}
	0 & 1\\
	1 & 0\end{array}\right).
\end{align} 
From Figure~\ref{fig:rej_rates_periodic}, we see that the belief converges to zero but beliefs oscillate a lot more about the mean rate of rejection as compared to the case of an aperiodic network given in equation (\ref{eq:aperiodic_w}).

Even though nodes do not have a positive self-weight ($W_{ii}$), the new information (through observations) entering at every node reaches its neighbors and gets dispersed in throughout the network; eventually reaches the node. Hence, nodes learn even when the network is periodic as long as it remains strongly connected. 

\subsubsection{Eigenvector Centrality and Extent of distinguishability}
From Theorem~\ref{thm:expconv}, we know that a larger weighted sum of the KL divergences, \emph{i.e.} a larger network divergence, $K(\theta_M, \theta_k)$, yields a better rate of rejecting hypothesis $\theta_k$. We look at a numerical example to show this. 

\begin{example}
Let $\Theta = \{ \theta_1, \theta_2, \theta_3, \theta_4, \theta_5\}$ and $\theta^* = \theta_4$. Consider a set of 25 nodes which are arranged in $5 \times 5$ array to form a grid. We obtain a grid network by connecting every node to its adjacent nodes. We define the weight matrix as,
\begin{align}
\label{eq:25grid}
	W_{ij} = \left\{ \begin{array}{cc}
	\frac{1}{\left| \mathcal{N}(i) \right|}, &\text{if } j \in \mathcal{N}(i)\\
	0, &\text{otherwise} \end{array}\right.
\end{align}
Consider an extreme scenario where only one node can distinguish true hypothesis $\theta_1$ from the rest and to the remaining nodes in the network all hypotheses are observationally equivalent \emph{i.e.} $\bar{\Theta}_i = \Theta$ for $24$ nodes and $\bar{\Theta}_i = \{ \theta_1 \}$ for only one node. We call that one node which can distinguish the true hypothesis from other hypotheses as the ``informed node'' and the rest of the nodes called the ``non-informed nodes''. 
\end{example}

\noindent 
For the weight matrix in equation (\ref{eq:25grid}), the eigenvector centrality of node $i$ is proportional to $\mathcal{N}(i)$, which means in this case, more number of neighbors implies higher social influence. This implies that the corner nodes namely node 1, node 5, node 20 and node 25 at the four corners of the grid have least eigenvector centrality among all nodes. Hence, they are least influential. The nodes on four edges have a greater influence than the corner nodes. Most influential nodes are the ones with four connections, such as node 13 which is located in third row and third column of the grid. It is also the central location of the grid. 

\begin{figure} [h!] 
  \centering
    \includegraphics[width=0.45\textwidth]{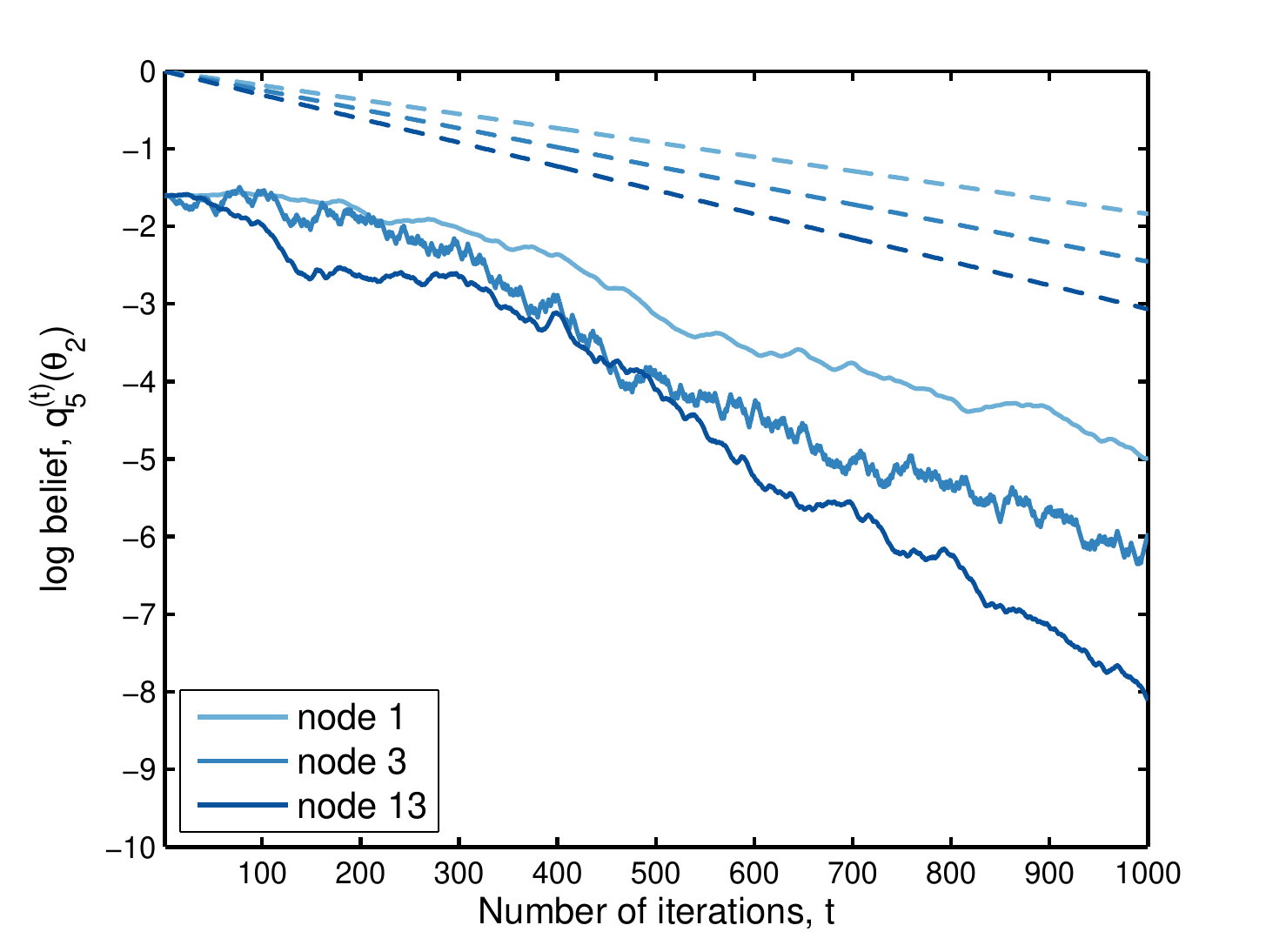}
    \caption{Figure illustrates the manner in which rate of rejection of $\theta_2$ at node 5 is influenced by varying the location of an informed node. As seen here when the informed node is more central \emph{i.e.} at node 13, rate of rejection is fastest and when the informed node is at the corner node 1, rate of rejection is slowest.}
    \label{fig:grid_informed_node_location}
\end{figure}

Figure~\ref{fig:grid_informed_node_location} shows the variation in the rate of rejection of $\theta_2$ of node 5 as the location of informed node changes. We see that if the informed node is at the center of the grid then the rate of rejection is fastest and the rate is slowest when the informed node is placed at a corner. In other words, rate of convergence is highest when the most influential node in the network has high distinguishability. 

\subsection{Factors influencing Concentration}
\label{subsec:conc_results}

\begin{figure} [h!] 
  \centering
    \includegraphics[width=0.45\textwidth]{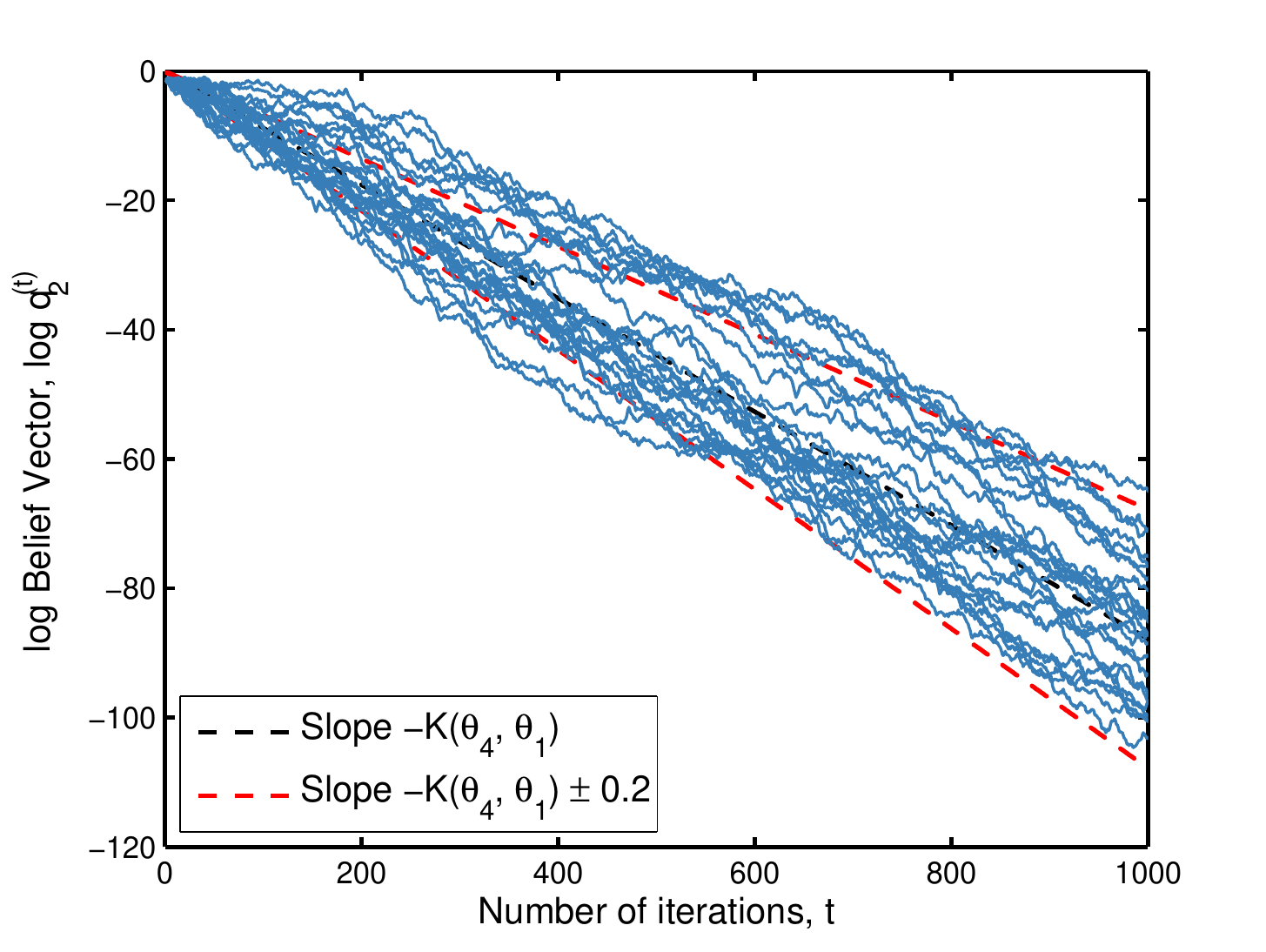}
    \caption{Figure shows the decay of belief of $\theta_1$ (wrong hypothesis) of node 2 for 25 instances. We see that the number of sample paths on which the rate of rejecting $\theta_1$ deviates more than $\eta = 0.1$ reduces as the number of iterations increase.}
    \label{fig:LDP_graph1}
\end{figure}

\begin{figure} [h!] 
  \centering
    \includegraphics[width=0.45\textwidth]{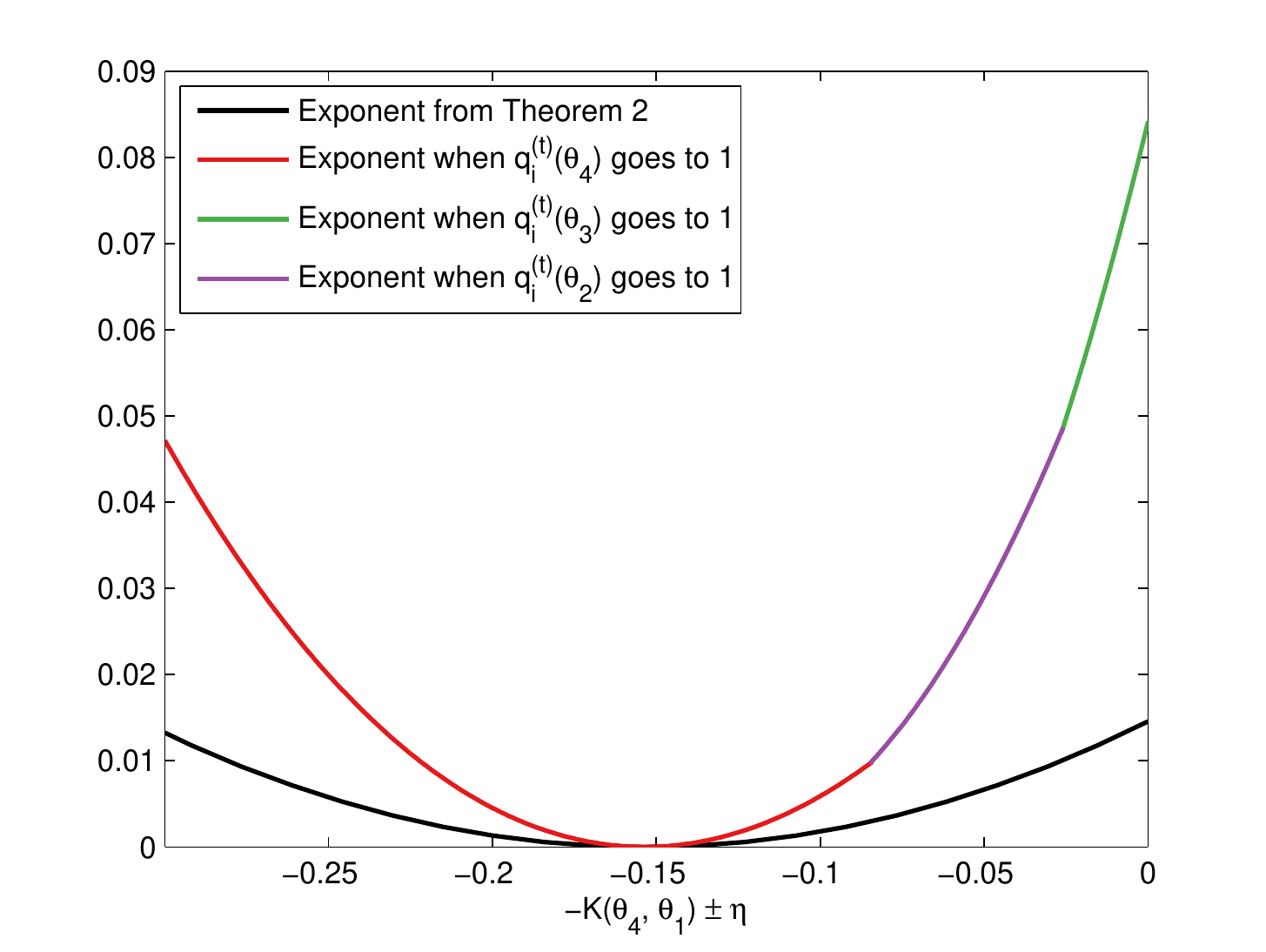}
    \caption{\add{Figure shows the asymptotic exponent with which the probability of events where rate of rejecting $\theta_1$ deviates by $\eta$ from $K(\theta_4, \theta_1)$; $\theta_4$ is the true hypothesis. The black curve shows the asymptotic exponent as characterized by Theorem~\ref{thm:expconvhoeffding}. The colored curve shows the exact asymptotic exponent as characterized by Theorem~\ref{thm:ldp_vector_belief}, where the exponent depends on the hypothesis to which the learning rule is converging. This shows that small deviations from $K(\theta_4, \theta_1)$ occur when the learning rule is converging to $\theta_4$ and larger deviations occur when the learning rule is converging to a wrong hypothesis. }}
    \label{fig:ldp_exponents}
\end{figure}

Now to examine the results from Theorem~\ref{thm:expconvhoeffding} and Theorem~\ref{thm:ldp_vector_belief}, we go back to Example~\ref{ex:2node}, where two nodes are in a strongly connected aperiodic network given by equation (\ref{eq:aperiodic_w}). Observation model for each node is defined as follows. For node 1, $\dist{1}{\cdot}{\theta_1} = \dist{1}{\cdot}{\theta_3}  \sim \Ber(\frac{4}{5})$ and $\dist{1}{\cdot}{\theta_2} = \dist{1}{\cdot}{\theta_4}  \sim \Ber(\frac{1}{4})$, and for node 2,  $\dist{2}{\cdot}{\theta_1} = \dist{2}{\cdot}{\theta_2} \sim \Ber(\frac{1}{3})$ and $\dist{2}{\cdot}{\theta_3} = \dist{2}{\cdot}{\theta_4}  \sim \Ber(\frac{1}{4})$. \add{Figure~\ref{fig:LDP_graph1} shows the exponential decay of $\theta_1$ for 25 instances. We see that the number of sample paths that deviate more than $\epsilon = 0.1$ from $K(\theta_4, \theta_1)$ decrease with number of iterations. Theorem~\ref{thm:expconvhoeffding} characterizes the asymptotic rate at which the probability of such sample paths vanishes when the log-likelihoods are bounded. This asymptotic rate is given as a function of $L$ and period of the network. From Corollary~\ref{coll:LDPtoHoeffding} we have that the rate given by Theorem~\ref{thm:expconvhoeffding} is loose for aperiodic networks. A tighter bound which utilizes the complete observation model is given by Theorem~\ref{thm:ldp_vector_belief}. Figure~\ref{fig:ldp_exponents} shows the gap between the rates.
}

\add{
Figure~\ref{fig:ldp_exponents} shows the rate at which the probability of sample paths deviating from rate of rejection can be thought of as operating in three different regimes. Here, each regime is denotes to the hypothesis to which the learning rule is converging. In order to see this consider the rate function of $\theta_1$, i.e.  $J_1(\cdot)$ from Corollary~\ref{coll:ldp_single_belief},
\begin{align*}
J_1(\mbf y) = \inf_{x \in \mathbb{R}^3: g(\mbf x) = y} I(\mbf x), \forall y \in \mathbb{R}.
\end{align*}
Behavior of the rate function $J_1(\cdot)$ depends on the function $g_1(\mbf x) = x_1 - \max\{0, x_1, x_2, x_3\}$. Whenever $g_1(\mbf x) = x_1$, the rate function is $I_1(\cdot)$. This shows that whenever there is a deviation of $x - k(\theta_4, \theta_1)$ from the rate of rejection of $\theta_1$, the sample paths that vanish with slowest exponents are those for which $\frac{1}{t} \log \frac{\est{i}{t}(\theta_1)}{\est{i}{t}(\theta_4)} < 0$ as $t \to \infty$. In other words, small deviations occur when the learning rule is converging to true hypothesis $\theta_4$ and they depend on $I_1(\cdot)$ (and hence $\theta_1$) alone. Whereas large deviations occur when the learning rule is mistakenly converging to a wrong hypothesis and hence, the rate function depends on $\theta_1$ and the wrong hypothesis to which the learning rule is converging. 
}

\add{
\subsection{Learning with Communication Constraints}
\label{sec:prac_issues}

Now, we consider a variant of our learning rule where the communication between the nodes is quantized to belong to a predefined finite set. Each node $i$ starts with an initial \add{private belief vector} $\mathbf{\est{i}{0}}$ and at each time $t = 1, 2,\ldots$ the following events happen:
	\begin{enumerate}
	\item Each node $i$ draws a conditionally i.i.d observation $\samp{i}{t} \sim \dist{i}{\cdot}{\theta_M}$.
	\item Each node $i$ performs a local Bayesian update on $\mathbf{\est{i}{t-1}}$ to form $\mathbf{\bel{i}{t}}$ using the following rule.
	 For each $k \in [M]$, 
		\begin{align}
		\bel{i}{t}(\theta_k) = \frac{ \dist{i}{\samp{i}{t}}{\theta_k} \est{i}{t-1}(\theta_k) }{ \sum_{a \in [M]} \dist{i}{\samp{i}{t}}{\theta_a} \est{i}{t-1}(\theta_a) }.
		\end{align}
	\item Each node $i$ sends the message $\msg{i}{t}(\theta_k) =  \left[ D \bel{i}{t} (\theta_k) \right]$, for all $k \in [M]$, to all nodes $j$ for which $i \in \mathcal{N}(j)$, where $D \in \mathbb{Z}^+$ and 
	\begin{align}
	\label{eq:quant_bel}
	[x] = 
	\left \{
	\begin{array}{ll}
	\lfloor x \rfloor + 1, & \text{ if } x > \lfloor x \rfloor +0.5,\\
	\lfloor x \rfloor, & \text{ if } x \leq \lfloor x \rfloor + 0.5,
	\end{array}
	\right.
	\end{align}
where $\lfloor x \rfloor$ denotes the largest integer less than $x$.
	\item Each node $i$ normalizes the beliefs received from the neighbors $\mathcal{N}(i)$ as 
	\begin{align}
	\tilde{Y}^{(t)}_i(\theta_k) = 
	\frac{\msg{i}{t}(\theta_k)}{\sum_{a \in [M]}\msg{i}{t}(\theta_a)},
	\end{align}	 
	and updates its private belief of $\theta_k$, for each  $k \in [M]$,
		\begin{align}
		\est{i}{t}(\theta_k) = \frac{ \exp \left( \sum_{j = 1}^{n} W_{ij} \log \tilde{Y}^{(t)}_i(\theta_k) \right)
			}{
			\sum_{a \in [M]} \exp \left( \sum_{j = 1}^{n} W_{ij} \tilde{Y}^{(t)}_i(\theta_a) \right)
			}.
		\end{align}
	\end{enumerate}
In the above learning rule, the belief on each  hypothesis belongs to a set of size $D+1$. Hence transmitting the entire belief vector, i.e., transmitting the entire message requires $M\log(D+1)$ bits.   	

Note that all of our simulations so far, we have used MATLAB. This means that in our simulations, we have relied on the default of  64-bit precision of the commonly used double-precision binary floating-point format to represent the belief on each hypothesis. This means our simulations can be interpreted as limiting the communication links to support 64 bits, or equivalently 8 bytes, per hypothesis 
per unit of time. Our previous experimental results show a close match between experiment and analysis using this level of quantization. In the next example we show the impact of coarser quantization in the following example.

\begin{example}
\label{ex:guassian_sensornet}
Consider a network with low cost radar or ultrasound sensors whose aim is to find the location of a target. Each sensor can sense the target's location along one dimension only, whereas the target location is a point in three-dimensional space. Consider the configuration in Figure~\ref{fig:sensor_network}: there are two nodes along each of the three coordinate axes at locations $[\pm 2, 0, 0]$, $[0,\pm 2, 0]$, and $[0,0,\pm 2]$. The communication links are given by the directed arrows. Nodes located on the x-axis can sense whether x-coordinate of the target lies in the interval $(-2,-1]$ or in the interval $(-1,0)$ or in the interval $[0, 1)$ or in the interval $[1,2)$. If a target is located in the interval $(-\infty, -2]\cup[2, \infty)$ on the x-axis then no node can detect it. Similarly nodes on y-axis and z-axis can each distinguish between 4 distinct non-intersecting intervals on the y-axis and the z-axis respectively. Therefore, the total number of hypotheses is $M = 4^3 = 64$.

The sensors receive signals which are three dimensional Gaussian vectors whose mean is altered in the presence of a target. In the absence of a target, the ambient signals have a Gaussian distribution with mean $[0,0,0]$. For the sensor node along $x$-axis located at $[2,0,0]$, if the target has $x$-coordinate $\theta_x \in (-2,2)$, the mean of the sensor's observation is $[ \lfloor 3 + \theta_x \rfloor, 0, 0 ]$. If a target is located in $(-\infty, -2]\cup[2, \infty)$ on the x-axis, then the mean of the Gaussian observations is $[0, 0, 0]$. Local marginals of the nodes along y-axis and z-axis are described similarly, i.e., as the target moves away from the node by one unit the signal mean strength goes by one unit. For targets located at a distance four units and beyond the sensor cannot detect the target. In this example, suppose $\theta_1$ is the true hypothesis.
\end{example}

\begin{figure} [h!] 
  \centering
    \includegraphics[width=0.45\textwidth]{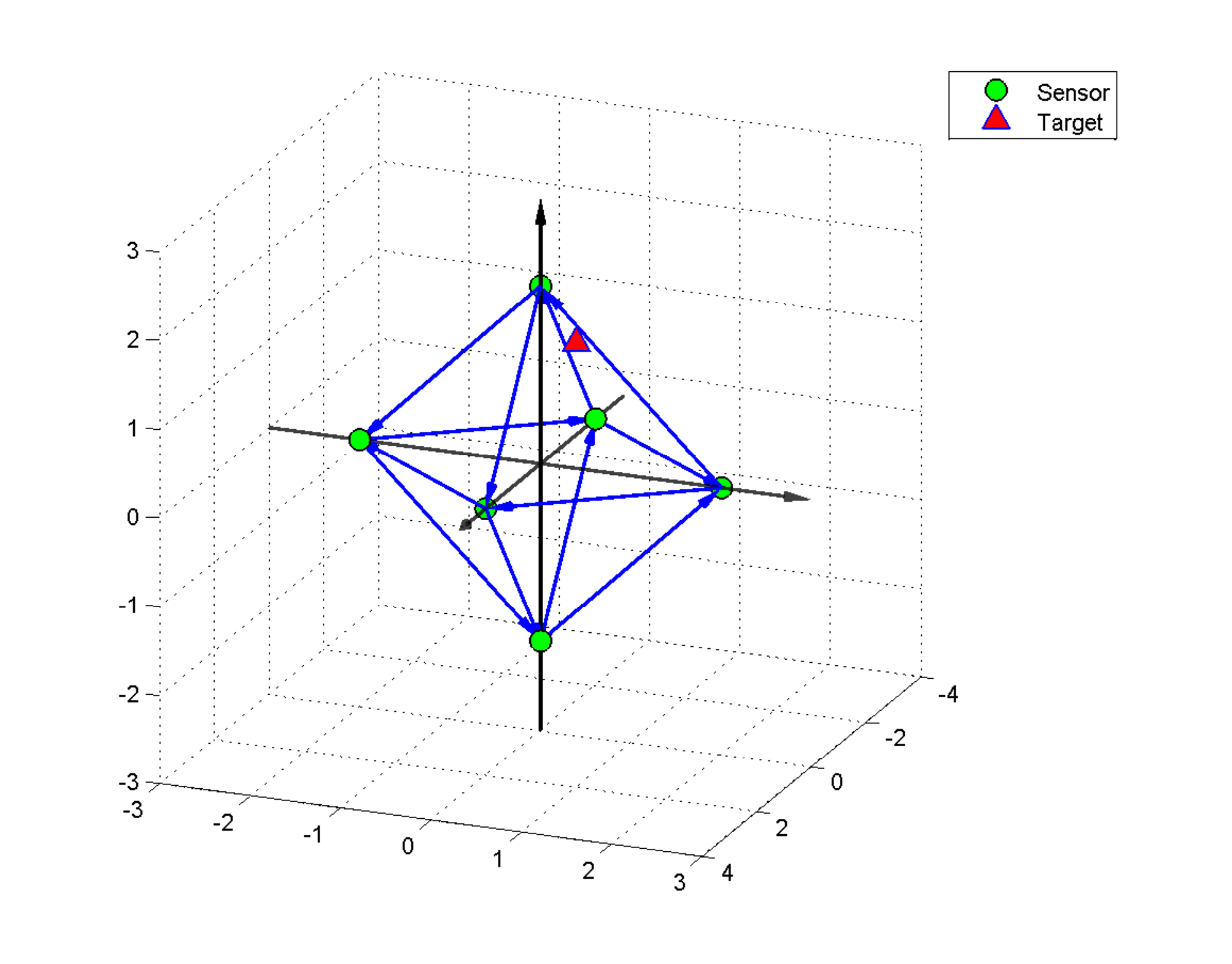}
    \caption{\add{Figure shows a sensor network where each node is a low cost radar that can sense along the axis it is placed and not the other. The directed edges indicate the directed communication between the nodes. Through cooperative effort the nodes aim to learn location of the target in 3-dimensions.}}
    \label{fig:sensor_network}
\end{figure}

\begin{figure} [h!] 
  \centering
    \includegraphics[width=0.45\textwidth]{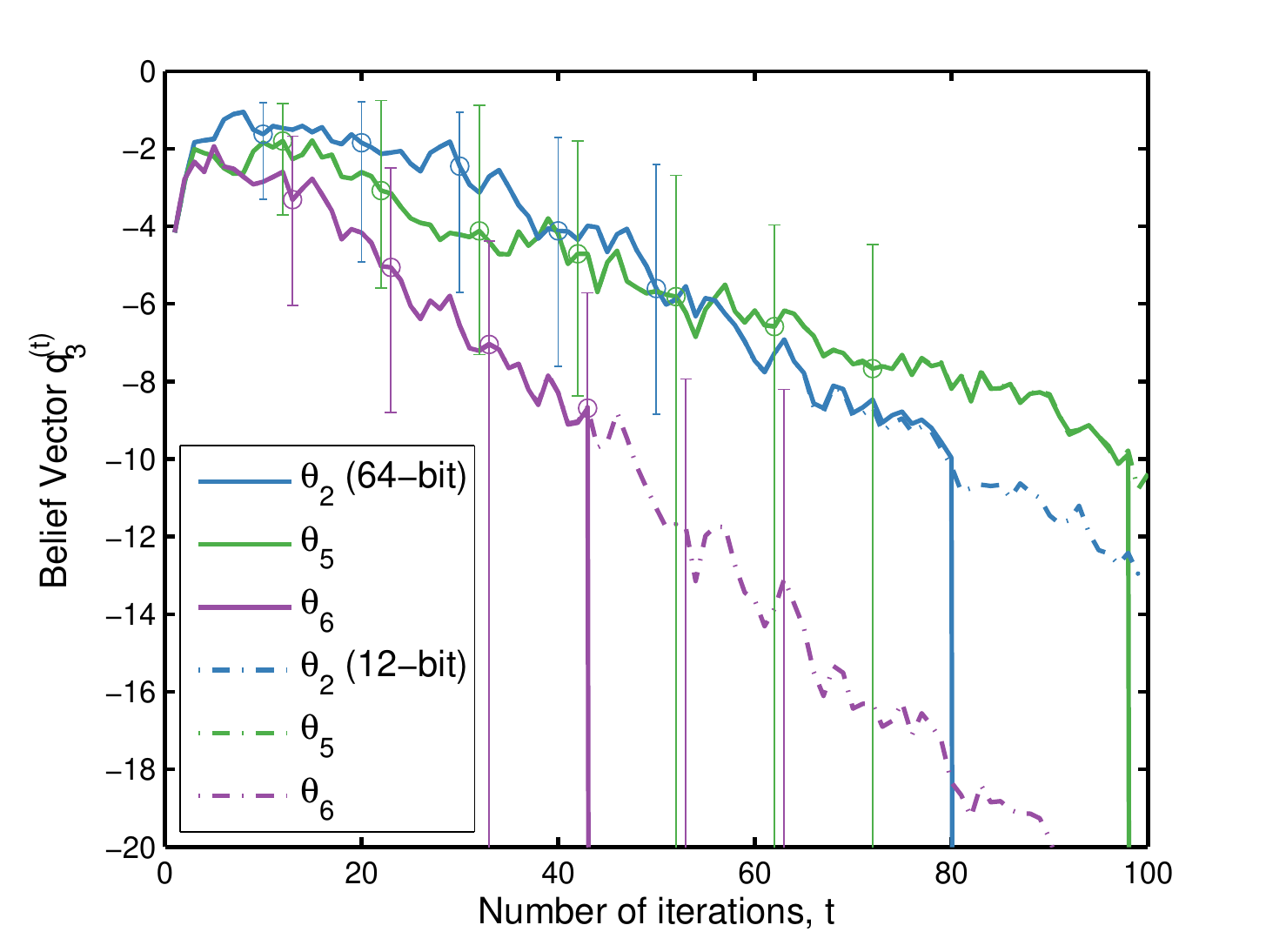}
    \caption{\add{The solid lines in figure show the evolution of the log beliefs of node 3 with time for hypotheses $\theta_2$, $\theta_5$ and $\theta_6$ when links support a maximum of 12 bits per hypothesis per unit time. This is compared with the evolution of the log beliefs with no rate restriction case (dotted lines) which translates a maximum of 64 bits per hypothesis per unit time. Figure also shows the confidence intervals around log beliefs over 500 instances of learning rule with 12 bits per hypothesis. We see the learning rule with link rate 12 bits per hypothesis converges in all the instances.}}
    \label{fig:12_bit_rate}
\end{figure}

\begin{figure} [h!] 
  \centering
    \includegraphics[width=0.45\textwidth]{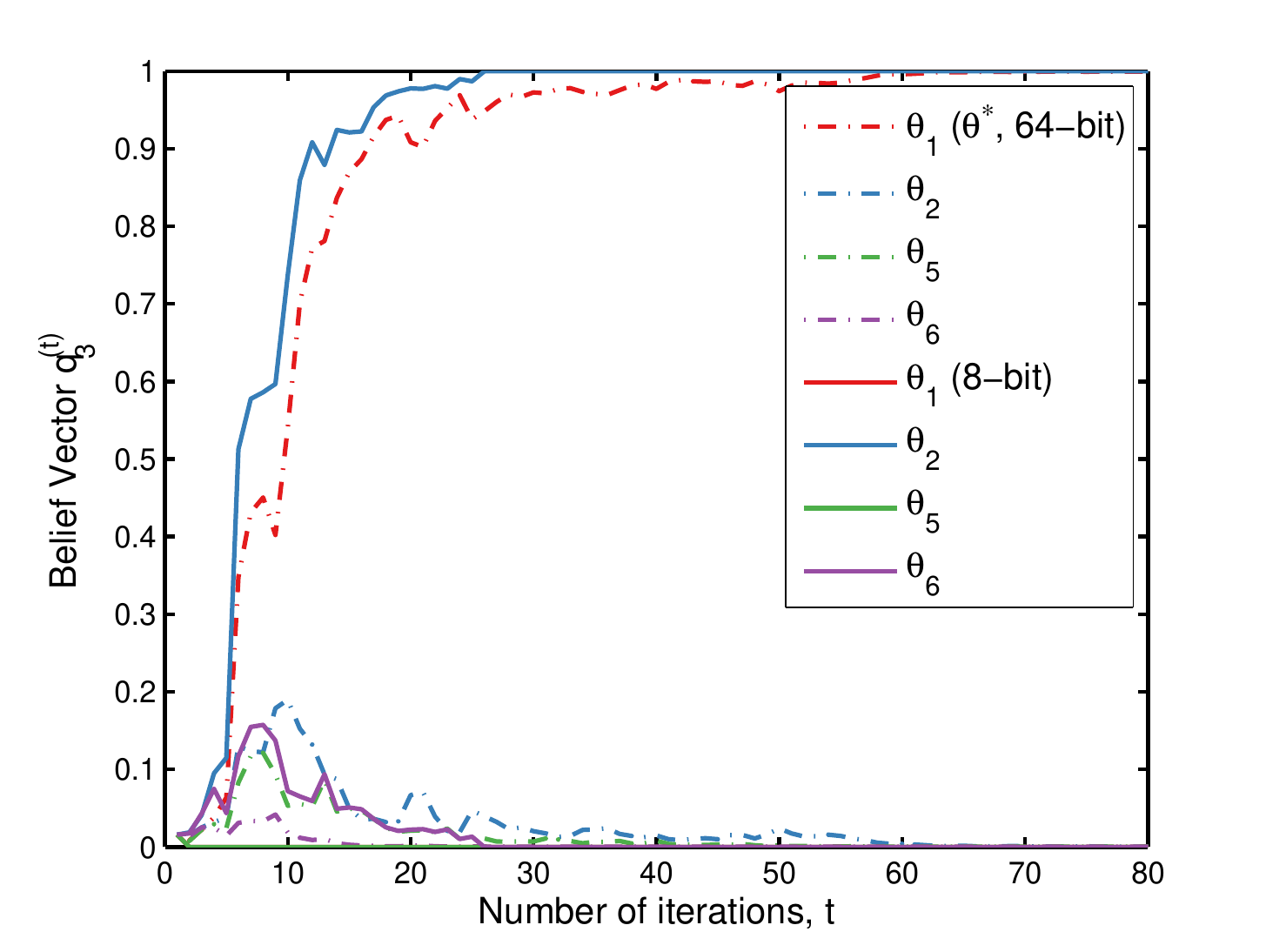}
    \caption{\add{The solid lines in the figure show the evolution of the log beliefs of node 3 with time for hypotheses $\theta_2$, $\theta_5$ and $\theta_6$ when links support a maximum of 8 bits per hypothesis per unit time. This is compared with the evolution of the log beliefs with no rate restriction case (dotted lines) which translates a maximum of 64 bits per hypothesis per unit time. For this sample path, we see that learning rule converges to a wrong hypothesis $\theta_5$ when the communication is restricted to 8 bits per hypothesis.}}
    \label{fig:8_bit_wrong}
\end{figure}

\begin{figure} [h!] 
  \centering
    \includegraphics[width=0.45\textwidth]{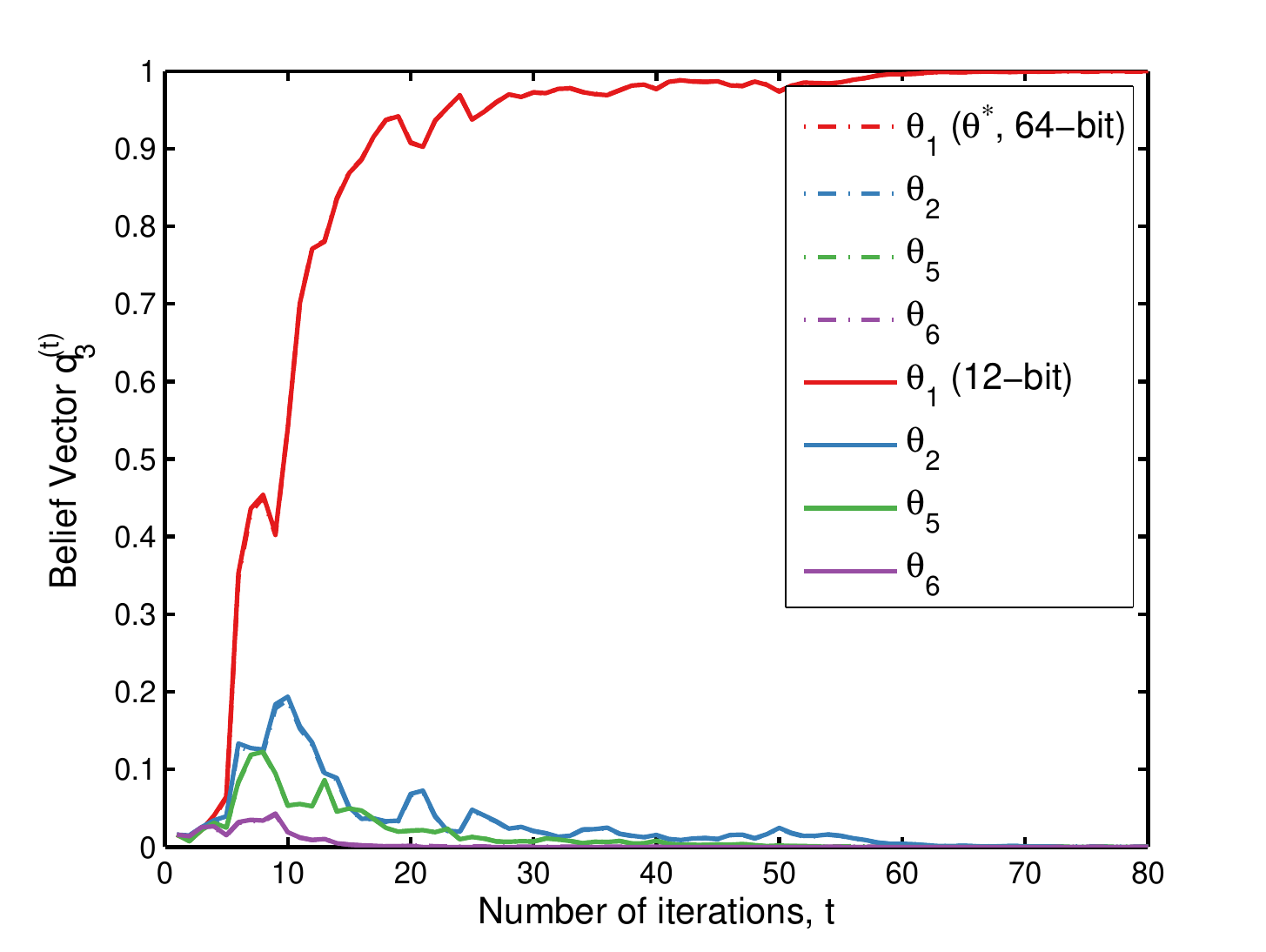}
    \caption{\add{The solid lines in figure show the evolution of the beliefs of node 3 with time for hypotheses $\theta_2$, $\theta_5$ and $\theta_6$ when links support a maximum of 12 bits per hypothesis per unit time. This is compared with the evolution of the beliefs with no rate restriction case (dotted lines) which in our simulations translates to the case when the links support a maximum of 64 bits per hypothesis per unit time. On the same sample path in Figure~\ref{fig:8_bit_wrong}, we see that learning rule converges to true hypothesis when the communication is restricted to 12 bits per hypothesis.}}
    \label{fig:12_bit}
\end{figure}


Consider $D = 2^{12}-1$ which implies that belief on each hypothesis is of size 12 bits or equivalently 1.5 bytes. Figure~\ref{fig:12_bit_rate} shows evolution of log beliefs of node 3 for hypotheses for $\theta_2$, $\theta_5$ and $\theta_6$ for 500 instances when the link rate is limited to 1.5 bytes per hypothesis per unit time. We see that the learning rule converges to the true hypotheses on all 500 instances. Now, consider $D = 2^{8}-1$ which implies that belief on each hypothesis is of size 8 bits or equivalently 1 byte. Figure~\ref{fig:8_bit_wrong} shows the evolution of beliefs of node 3 for hypotheses $\theta_2$, $\theta_5$ and $\theta_6$ when the link rate is limited to 1 byte per hypothesis per unit time. We see that the learning rule converges to a wrong hypothesis $\theta_2$. Whereas, on the same sample path in Figure~\ref{fig:12_bit} we see that if the link rate is 1.5 bytes per hypothesis per unit time, the learning rule converges to true hypothesis. This happens because on every sample path our learning rule has an initial transient phase where beliefs may have large fluctuations during which the belief on true hypothesis may get close to zero. For low link rates, the value of $D$ is small and when the belief on true hypothesis even though strictly positive becomes less than $\frac{1}{2D}$, it gets quantized to zero. Recall from Assumption~\ref{assume:initial_est}, our learning rule when a belief goes to zero, propagates the zero belief to all subsequent time instants. This shows that as we increase the value of $D$, i.e., as we increase link rate, the quantized learning rule is more robust to the initial fluctuations but there are certain samples on which it may converge to a wrong hypothesis.

In our simulations overall, we observe that for both Examples~\ref{ex:2node} and~\ref{ex:guassian_sensornet}, when link rates are greater than or equal to 1.5 bytes per hypothesis per unit time the learning rule converges for all instances and its performance coincides with the prediction of our the analysis under the assumption of perfect links.}


\section{Conclusion} \label{sec:discussion}

\add{
In this paper we study protocols in which a network of nodes make observations and communicate in order to collectively learn an unknown fixed global hypothesis that statistically governs the distribution of their observations. Our learning rule performs local Bayesian updating followed by averaging log-beliefs. We show that our protocol guarantees exponentially fast convergence to the true hypothesis with probability one. 
We showed the rate of rejection of any wrong hypothesis has an explicit characterization in terms of the local divergences and network topology. Furthermore, under the (mild) Assumption~\ref{assume:finite_lmgf}, we provide an asymptotically tight characterization of rate of concentration for the rate of rejection. This assumption admits a broad class of distributions with unbounded support such as Gaussian mixtures. 

Future work should consider more practical limitations in our setting. Our experimental results indicate that if the rate limit on link is sufficiently high, our protocol will still be successful. This indicates our learning rule is a first step towards a more realistic study of distributed hypothesis testing with more practical constraints on communication. An open question is the minimum data rate on communication link to ensure convergence of learning to true hypothesis: an analytic study of the learning rule at low data rates can make a connection to previous results on rate-limited distributed hypothesis testing. 
}

\appendix
\subsection{Proof of Theorem~\ref{thm:expconv}}
We begin with the following recursion for each node $i$ and $k \in [M-1]$;
	\begin{align}
	&\log \frac{ \est{i}{t}(\theta_M) }{ \est{i}{t}(\theta_k) } 
	\nonumber
	\\
	&=\sum_{j = 1}^{n} W_{ij} \log \frac{ \bel{j}{t}(\theta_M) }{ \bel{j}{t}(\theta_k) } \nonumber
	\\
	& 
	= \sum_{j = 1}^{n} W_{ij} \left(
		  \log \frac{ \dist{j}{\samp{j}{t}}{\theta_M} }{ \dist{j}{\samp{j}{t}}{\theta_k} }
		+
		\log \frac{ \est{j}{t-1}(\theta_M) }{ \est{j}{t-1}(\theta_k) }
		\right). 
	\label{eq:recursion}
	\end{align} 
where the first and the second equalities follow from~\eqref{eq:estimate} and~\eqref{eq:bayes}, respectively. 
Now for each node $j$ we rewrite $\log \frac{\est{j}{\cdot}(\theta_M)}{\est{j}{\cdot}(\theta_k)}$  in terms of node $j$'s neighbors and their samples at the previous instants. We can expand in this way until we express everything in terms of the samples collected and the initial estimates. 
Noting that $W^{t}(i,j) = \sum_{i_{t-1} = 1}^{n} \ldots \sum_{i_1 = 1}^{n} W_{ii_{1}} \ldots W_{i_{t - 1}j}$, it is easy to check that equation (\ref{eq:recursion}) can be further expanded as
 \begin{align}
&\log \frac{\est{i}{t}(\theta_M)}{\est{i}{t}(\theta_k)} \nonumber
	\\
 	& =   \sum_{j = 1}^{n}\sum_{\tau = 1}^{t} W^{\tau}(i, j) \log \frac{ \dist{j}{\samp{j}{t - \tau + 1}}{\theta_M} }{ \dist{j}{\samp{j}{t - \tau + 1}}{\theta_k}}\nonumber
 	\\
 	&\hspace{0.5cm}+  \sum_{j = 1}^{n}W^{t}(i, j) \log \frac{\est{j}{0}(\theta_M) }{ \est{j}{0}(\theta_k)}.
\label{equation:compactform}
 \end{align}
 Now divide by $t$ and take limit as $t \to \infty$
\begin{align}
 &\lim_{t \to \infty} \frac{1}{t}\log \frac{\est{i}{t}(\theta_M)}{\est{i}{t}(\theta_k)} \nonumber
 	\\
 	&= \lim_{t \to \infty} \frac{1}{t}\sum_{j = 1}^{n}\sum_{\tau = 1}^{t} W^{\tau}(i, j) \log \frac{ \dist{j}{\samp{j}{t - \tau + 1}}{\theta_M} }{ \dist{j}{\samp{j}{t - \tau + 1}}{\theta_k}}				\nonumber
 	\\
 	& \hspace{0.5cm}+ \lim_{t \to \infty} \frac{1}{t}\sum_{j = 1}^{n}W^{t}(i, j) \log \frac{\est{j}{0}(\theta_M) }{ \est{j}{0}(\theta_k)}.
\label{eq:compactform_lim}
 \end{align}
From Assumption~\ref{assume:initial_est}, the prior $\est{j}{0}(\theta_k)$ is strictly positive for every node $j$ and every $k \in [M]$. Since $W^{t}(i,j) \leq 1$, we have  
\begin{equation}
\lim_{t \to \infty}\frac{1}{t} \left\{\sum_{j = 1}^{n}W^{t}(i, j) \log \frac{ \est{j}{0}(\theta_M)}{ \est{j}{0}(\theta_k)} \right\}  = 0.
\label{eq:lim_prior}
\end{equation}  

Let $W$ be periodic with period $d$. If $W$ is aperiodic, then the same proof still holds by putting $d = 1$. Now, we fix node $i$ as a reference node and for every $r \in [d]$, define
\begin{equation*}
A_r = \{j \in [n]: W^{md + r}(i,j) > 0 \text{ for some } m \in \mathbb{N}\}.
\end{equation*}
 In particular, $(A_1, A_2, \ldots, A_{d})$ is a partition of $[n]$; these sets form cyclic classes of the Markov chain. Fact~\ref{fact:stationarydist} implies that for every $\delta > 0$, there exists an integer $N$ which is function of $\delta$ alone, such that for all $m \geq N$, for some fixed $r \in [d]$ if $ j \in A_r$, then
\begin{equation}
\left| W^{md + r}(i, j) - v_j d \right| \leq \delta
\label{eq:stationarydist1}
\end{equation}
and if $j \not \in A_r$
\begin{equation}
0 \leq W^{md + r}(i, j) \leq \delta.
\label{eq:stationarydist2}
\end{equation}
Using this the first term in equation (\ref{eq:compactform_lim}) can be decomposed as follows
\begin{align}
&\lim_{t \to \infty} \frac{1}{t}\sum_{j = 1}^{n}\sum_{\tau = 1}^{t} W^{\tau}(i, j) \log \frac{ \dist{j}{\samp{j}{t - \tau + 1}}{\theta_M} }{ \dist{j}{\samp{j}{t - \tau + 1}}{\theta_k}}	 \nonumber
	\\		
	&= \lim_{t \to \infty} \frac{1}{t}\sum_{j = 1}^{n}\sum_{\tau = 1}^{Nd - 1} W^{\tau}(i, j) \log \frac{ \dist{j}{\samp{j}{t - \tau + 1}}{\theta_M} }{ \dist{j}{\samp{j}{t - \tau + 1}}{\theta_k}}	 \nonumber
	\\
	& \hspace{0.5cm} + \lim_{t \to \infty} \frac{1}{t}\sum_{j = 1}^{n}\sum_{\tau = Nd}^{t} W^{\tau}(i, j) \log \frac{ \dist{j}{\samp{j}{t - \tau + 1}}{\theta_M} }{ \dist{j}{\samp{j}{t - \tau + 1}}{\theta_k}}.	
\label{eq:sum1}	
 \end{align}

\noindent
Using triangle inequality and the fact that $W^{\tau}(i,j) \leq 1$ for every $\tau \in \mathbb{N}$ we have
\begin{align*}
&\left| \lim_{t \to \infty} \frac{1}{t}\sum_{\tau = 1}^{Nd - 1}W^{\tau}(i,j)\log \frac{ \dist{j}{\samp{j}{t - \tau}}{\theta_M} }{ \dist{j}{\samp{j}{t - \tau}}{\theta_k}} \right| \nonumber
	\\
 	&\leq 
	\lim_{t \to \infty} \frac{1}{t}\sum_{\tau = 1}^{Nd - 1} \left| \log \frac{ \dist{j}{\samp{j}{t - \tau}}{\theta_M} }{ \dist{j}{\samp{j}{t - \tau}}{\theta_k}} \right|. 
\end{align*}
For every $j \in [n]$, $\log \frac{ \dist{j}{X_j}{\theta_M} }{ \dist{j}{X_j}{\theta_k}}$
is integrable, implying $\left| \log\frac{ \dist{j}{X_j}{\theta_M} }{ \dist{j}{X_j}{\theta_k}} \right|$ is almost surely finite. This implies that
\begin{align}
\lim_{t \to \infty} \frac{1}{t}\sum_{\tau = 1}^{Nd - 1}W^{\tau}(i,j)\log \frac{ \dist{j}{\samp{j}{t - \tau}}{\theta_M} }{ \dist{j}{\samp{j}{t - \tau}}{\theta_k}} = 0 \quad \P \mbox{-a.s.}
\label{eq:finite_as}
\end{align}
 
\noindent
Using (\ref{eq:lim_prior}) and (\ref{eq:finite_as}), equation (\ref{eq:sum1}) becomes
\begin{align*}
&\lim_{t \to \infty} \frac{1}{t}\log \frac{\est{i}{t}(\theta_M)}{\est{i}{t}(\theta_k)} \nonumber
 	\\
 	&= \lim_{t \to \infty}  \frac{1}{t}\sum_{j = 1}^{n}\sum_{\tau = Nd}^{t} W^{\tau}(i, j) \log \frac{ \dist{j}{\samp{j}{t - \tau + 1}}{\theta_M} }{ \dist{j}{\samp{j}{t - \tau + 1}}{\theta_k}},
\end{align*} with probability one. It is straightforward to see that the above equation can be rewritten as
\begin{align*}
&\lim_{t \to \infty} \frac{1}{t}\log \frac{\est{i}{t}(\theta_M)}{\est{i}{t}(\theta_k)} \nonumber
	\\
	& = \lim_{T \to \infty}  \frac{1}{Td}\sum_{j = 1}^{n}\sum_{m = N}^{T - 1}\left\{ \sum_{r = 1}^{d} W^{md + r}(i, j) \times \right.\nonumber
	\\
	&\hspace{3cm} \left. \log \frac{ \dist{j}{\samp{j}{Td - md - r+ 1}}{\theta_M} }{ \dist{j}{\samp{j}{Td - md - r+ 1}}{\theta_k}}  \right \},
\end{align*} with probability one. For every $\delta > 0$ and $N$ such that for all $m \in N$ equations (\ref{eq:stationarydist1}) and (\ref{eq:stationarydist2}) hold true, using Lemma~\ref{lemma:interval} we get that
\[
\lim_{t \to \infty} \frac{1}{t}\log \frac{\est{i}{t}(\theta_M)}{\est{i}{t}(\theta_k)}
\]
with probability one lies in the interval with end points
\[
  K(\theta_M, \theta_k) - \frac{  \delta}{d} \sum_{j = 1}^{n} \expe\left[  \left| \log\frac{ \dist{j}{X_j}{\theta_M} }{ \dist{j}{X_j}{\theta_k}} \right| \right]
\]
and
\[
  K(\theta_M, \theta_k) + \frac{ \delta}{d} \sum_{j = 1}^{n} \expe\left[  \left| \log\frac{ \dist{j}{X_j}{\theta_M} }{ \dist{j}{X_j}{\theta_k}} \right| \right].
\]
Since this holds for any $\delta > 0$, we have
 \begin{equation*}
 \lim_{t \to \infty}\frac{1}{t}\log \frac{\est{i}{t}(\theta_M)}{\est{i}{t}(\theta_k)} =    K(\theta_M, \theta_k) \quad \P\mbox{-a.s.}
 \end{equation*} 
Hence, with probability one, for every $\epsilon > 0$ there exists a time $T^{\prime}$ such that $\forall t \geq T^{\prime}$, $\forall k \in [M-1]$ we have
 \begin{equation*}
 \left| \frac{1}{t}\log \frac{\est{i}{t}(\theta_M)}{\est{i}{t}(\theta_k)} -     K(\theta_M, \theta_k) \right| \leq \epsilon,
 \end{equation*} 
which implies
 \begin{equation*}
\frac{1}{1 + \displaystyle \sum_{k \in [M-1]}e^{-   K(\theta_M, \theta_k)t +\epsilon t}}
	\leq 
	\est{i}{t}(\theta_M)
	\leq 1.
 \end{equation*}
Hence we have the assertion of the theorem.

\begin{lemma}
\label{lemma:interval}
For a given $\delta > 0$ and for some $N \in \mathbb{N}$ for which equations (\ref{eq:stationarydist1}) and (\ref{eq:stationarydist2}) hold true for all $m \geq N$, the following expression 
\begin{align*}
&\lim_{T \to \infty} \frac{1}{Td}\sum_{j = 1}^{n}\sum_{m = N}^{T - 1}\left\{ \sum_{r = 1}^{d} W^{md + r}(i, j) \times \right.\nonumber
	\\
	&\hspace{3.5cm} \left. \log \frac{ \dist{j}{\samp{j}{Td - md - r+ 1}}{\theta_M} }{ \dist{j}{\samp{j}{Td - md - r+ 1}}{\theta_k}}  \right \}
\end{align*}
with probability one lies in an interval with end points
\[K(\theta_M, \theta_k) - \frac{\delta}{d} \sum_{j = 1}^{n}\expe\left[ \left| \log\frac{ \dist{j}{X_j}{\theta_M} }{ \dist{j}{X_j}{\theta_k}} \right| \right], 
\] and
\[
K(\theta_M, \theta_k) + \frac{\delta}{d} \sum_{j = 1}^{n}\expe\left[ \left| \log\frac{ \dist{j}{X_j}{\theta_M} }{ \dist{j}{X_j}{\theta_k}} \right| \right].
\]
\end{lemma}

\begin{IEEEproof}
To the given expression we add and subtract $v_j d$ from $W^{md + r}(i, j)$ for all $j \in A_r$ and we get
\begin{align}
&\lim_{T \to \infty} \frac{1}{Td}\sum_{j = 1}^{n}\sum_{m = N}^{T - 1}\left\{ \sum_{r = 1}^{d} W^{md + r}(i, j) \times \right.\nonumber
	\\
	&\hspace{3.5cm} \left. \log \frac{ \dist{j}{\samp{j}{Td - md - r+ 1}}{\theta_M} }{ \dist{j}{\samp{j}{Td - md - r+ 1}}{\theta_k}}  \right \}\nonumber
	\\
	& = \sum_{r = 1}^{d} \sum_{j \not \in A_r} \left\{  \lim_{T \to \infty} \frac{1}{Td} \sum_{m = N}^{T - 1}W^{md + r}(i, j) \right. \times \nonumber
	\\
	&\hspace{3.5cm}\left.\log \frac{ \dist{j}{\samp{j}{Td - md - r+ 1}}{\theta_M} }{ \dist{j}{\samp{j}{Td - md - r+ 1}}{\theta_k}}  \right \} \nonumber
	\\
	&\hspace{0.4cm} + \sum_{r = 1}^{d} \sum_{j \in A_r} \left\{ \lim_{T \to \infty} \frac{1}{Td}\sum_{m = N}^{T - 1}\left( W^{md + r}(i, j)- v_j d \right) \right.  \times \nonumber
	\\
	&\hspace{3.5cm} \left. \log \frac{ \dist{j}{\samp{j}{Td - md - r+ 1}}{\theta_M} }{ \dist{j}{\samp{j}{Td - md - r+ 1}}{\theta_k}} \right \} \nonumber
	\\
	&\hspace{0.5cm} + \sum_{r = 1}^{d} \sum_{j \in A_r}\left\{ \lim_{T \to \infty} \frac{1}{Td}\sum_{m = N}^{T-1} v_j d \right. \times \nonumber
	\\
	&\hspace{3.5cm}\left.\log \frac{ \dist{j}{\samp{j}{Td - md - r+ 1}}{\theta_M} }{ \dist{j}{\samp{j}{Td - md - r+ 1}}{\theta_k}} \right \} .
	\label{eq:sum2}
\end{align}

\noindent
For each $r$ and some $j \in A_r$, using equation (\ref{eq:stationarydist1}) and strong law of large numbers we have
%
%
\begin{align*}
& \left| \lim_{T \to \infty} \frac{1}{Td} \left\{ \sum_{m = N}^{T - 1} \left( W^{md + r}(i, j)- v_j d \right) \times \right. \right. \nonumber
	\\
	 & \hspace{3cm}\left. \left.\log \frac{ \dist{j}{\samp{j}{Td - md - r+ 1}}{\theta_M} }{ \dist{j}{\samp{j}{Td - md - r+ 1}}{\theta_k}} \right\} \right| \nonumber
	\\
	& \leq
	\frac{\delta }{d}\expe\left[  \left| \log\frac{ \dist{j}{X_j}{\theta_M} }{ \dist{j}{X_j}{\theta_k}} \right| \right] \quad \P\mbox{-a.s.}
\end{align*}
Similarly for $j \not \in A_r$, using equation (\ref{eq:stationarydist2}) we have
\begin{align*}
&\left| \lim_{T \to \infty} \frac{1}{Td}\sum_{m = N}^{T - 1}  W^{md + r}(i, j) \right. \times 
	\\
	&\hspace{3cm} \left. \log \frac{ \dist{j}{\samp{j}{Td - md - r+ 1}}{\theta_M} }{ \dist{j}{\samp{j}{Td - md - r+ 1}}{\theta_k}} \right|
	\\
	&\leq
	\frac{\delta }{d}\expe\left[  \left| \log\frac{ \dist{j}{X_j}{\theta_M} }{ \dist{j}{X_j}{\theta_k}} \right| \right] \quad \P\mbox{-a.s.}
\end{align*}
Again, by the strong law of large numbers we have
\begin{align*}
&\sum_{r = 1}^{d} \sum_{j \in A_r} v_j  \left\{ \lim_{T \to \infty} \frac{1}{T}\sum_{m = N}^{T-1} \log \frac{ \dist{j}{\samp{j}{Td - md - r+ 1}}{\theta_M} }{ \dist{j}{\samp{j}{Td - md - r+ 1}}{\theta_k}} \right \} \nonumber
	\\
	&= \sum_{r = 1}^{d} \sum_{j \in A_r}v_j \expe\left[ \log\frac{ \dist{j}{X_j}{\theta_M} }{ \dist{j}{X_j}{\theta_k}} \right] \nonumber
	\\
	&=K(\theta_M, \theta_k) \quad \P \mbox{-a.s.}
\end{align*}
Now combining this with equation (\ref{eq:sum2}) we have the assertion of the lemma.

\end{IEEEproof}

\subsection{Proof of Theorem~\ref{thm:expconvhoeffding}}

Recall the following equation
\begin{align}
&\lim_{t \to \infty}  \frac{1}{t} \log \frac{\est{i}{t}(\theta_M)}{\est{i}{t}(\theta_k)} \nonumber
	\\		
	&= \lim_{t \to \infty}  \frac{1}{t}\sum_{j = 1}^{n}\sum_{\tau = 1}^{Nd - 1} W^{\tau}(i, j) \log \frac{ \dist{j}{\samp{j}{t - \tau + 1}}{\theta_M} }{ \dist{j}{\samp{j}{t - \tau + 1}}{\theta_k}}	 \nonumber
	\\
	& \hspace{0.5cm} + \lim_{t \to \infty}  \frac{1}{t}\sum_{j = 1}^{n}\sum_{\tau = Nd}^{t} W^{\tau}(i, j) \log \frac{ \dist{j}{\samp{j}{t - \tau + 1}}{\theta_M} }{ \dist{j}{\samp{j}{t - \tau + 1}}{\theta_k}},	
\label{eq:sum1_thm2}	
 \end{align}
where $N$ is such that for all $m \geq N, m \in \mathbb{N}$ equation~\ref{eq:stationarydist1} and~\ref{eq:stationarydist2} are satisfied. For any fixed $t$, using Assumption~\ref{assume:boundedMG}, the first term in the summation on the right hand side of equation~\ref{eq:sum1_thm2} can be bounded as
\begin{align*}
&\left| \frac{1}{t}\sum_{j = 1}^{n}\sum_{\tau = 1}^{Nd - 1} W^{\tau}(i, j) \log \frac{ \dist{j}{\samp{j}{t - \tau + 1}}{\theta_M} }{ \dist{j}{\samp{j}{t - \tau + 1}}{\theta_k}}\right|
	\leq \frac{nNdL}{t}.
\end{align*}
Also, the second term in the summation on the right hand side of equation~\ref{eq:sum1_thm2} can be bounded as
\begin{align*}
&\left| \frac{1}{t}\sum_{j = 1}^{n}\sum_{\tau = Nd}^{t} W^{\tau}(i, j) \log \frac{ \dist{j}{\samp{j}{t - \tau + 1}}{\theta_M} }{ \dist{j}{\samp{j}{t - \tau + 1}}{\theta_k}} \right.
	\\
	& - \left. \sum_{r = 1}^{d} \sum_{j \in A_r} \frac{v_j }{Td}\sum_{m = 0}^{T - 1} \log \frac{ \dist{j}{\samp{j}{Td - md - r+ 1}}{\theta_M} }{ \dist{j}{\samp{j}{Td - md - r+ 1}}{\theta_k}} \right|
	\\
	& \leq \delta \frac{1}{Td}\sum_{m = 0}^{T - 1} \left| \log \frac{ \dist{j}{\samp{j}{Td - md - r+ 1}}{\theta_M} }{ \dist{j}{\samp{j}{Td - md - r+ 1}}{\theta_k}} \right|.
\end{align*}
Using Assumption~\ref{assume:boundedMG} we have
\begin{align*}
\frac{1}{Td}\sum_{m = 0}^{T - 1} \left| \log \frac{ \dist{j}{\samp{j}{Td - md - r+ 1}}{\theta_M} }{ \dist{j}{\samp{j}{Td - md - r+ 1}}{\theta_k}} \right| \leq \frac{L}{d}.
\end{align*}
Therefore, we have 
\begin{align*}
&\left|  \frac{1}{t} \log \frac{\est{i}{t}(\theta_M)}{\est{i}{t}(\theta_k)} \right.
	\\
	& - \left. \sum_{r = 1}^{d} \sum_{j \in A_r}   \frac{  v_j}{Td}\sum_{m = 0}^{T - 1} \log \frac{ \dist{j}{\samp{j}{Td - md - r+ 1}}{\theta_M} }{ \dist{j}{\samp{j}{Td - md - r+ 1}}{\theta_k}} \right|
	\\
	& \leq \frac{\delta   n L}{d}.
\end{align*}

\noindent
Applying Hoeffding's inequality (Theorem~2 of~\cite{Hoeffding:JSTOR1963}), equation (\ref{eq:sum1_thm2}) for $t \geq Nd$, for every $0 < \epsilon \leq K(\theta_M, \theta_k)$ can be written as 
\begin{align*}
\frac{1}{t}\log \frac{\est{i}{t}(\theta_M)}{\est{i}{t}(\theta_k)}
\leq K(\theta_M, \theta_k) - \epsilon + o\left( \frac{1}{t}, \delta \right),
\end{align*}
with probability at most  $\exp\left(-\frac{\epsilon^2T}{2 L^2}\right)$ where $o\left( \frac{1}{t}, \delta \right) = \frac{\delta   n L}{d} + \frac{nNdL}{t}$. Similarly, for $0 < \epsilon \leq L - K(\theta_M, \theta_k)$ we have
\begin{align*}
\frac{1}{t}\log \frac{\est{i}{t}(\theta_M)}{\est{i}{t}(\theta_k)}
\geq K(\theta_M, \theta_k) + \epsilon + o\left( \frac{1}{t}, \delta \right),
\end{align*}
with probability at most  $\exp\left(-\frac{\epsilon^2T}{2 L^2}\right)$ and for $\epsilon > L - K(\theta_M, \theta_k)$ we have
\begin{align*}
\frac{1}{t}\log \frac{\est{i}{t}(\theta_M)}{\est{i}{t}(\theta_k)}
\geq K(\theta_M, \theta_k) + \epsilon + o\left( \frac{1}{t}, \delta \right),
\end{align*}
with probability 0. 
Now, taking limit and letting $\delta$ go to zero, for $0 < \epsilon \leq K(\theta_M, \theta_k)$ we have
\begin{align*}
\lim_{t \to \infty} \frac{1}{t}\log \P \left(  \rho_{i}^{(t)}(\theta_k) - \rho_{i}^{(t)}(\theta_M) \leq  K(\theta_M, \theta_k) - \epsilon \right) 
\\
\leq -\frac{\epsilon^2}{2   L^2 d},
\end{align*}  
for $0 < \epsilon \leq L - K(\theta_M, \theta_k)$ we have
\begin{align*}
\lim_{t \to \infty} \frac{1}{t}\log \P \left(  \rho_{i}^{(t)}(\theta_k) - \rho_{i}^{(t)}(\theta_M) \geq  K(\theta_M, \theta_k) + \epsilon \right) 
\\
\leq -\frac{\epsilon^2}{2 L^2 d},
\end{align*} 
and for $\epsilon > L- K(\theta_M, \theta_k)$ we have
\begin{align*}
\lim_{t \to \infty} \frac{1}{t}\log \P \left(  \rho_{i}^{(t)}(\theta_k) - \rho_{i}^{(t)}(\theta_M) \geq  K(\theta_M, \theta_k) + \epsilon \right) 
\\
= -\infty,
\end{align*}
Since $\est{i}{t}(\theta_M) \leq 1$, all the events $\omega$ which lie in the set $\{\omega : \rho_{i}^{(t)}(\theta_k) \leq   K(\theta_M, \theta_k) - \epsilon\}$ also lie in the set $\{\omega : \rho_{i}^{(t)}(\theta_k)  \leq   K(\theta_M, \theta_k) - \epsilon + \rho_{i}^{(t)}(\theta_M)\}$. Hence, for every $0 < \epsilon \leq K(\theta_M, \theta_k)$ we have
\begin{align}
\label{eq:concentration_below_mean}
\lim_{t \to \infty} \frac{1}{t}\log \P \left(\rho_{i}^{(t)}(\theta_k) \leq    K(\theta_M, \theta_k) - \epsilon \right)
\leq 
-\frac{\epsilon^2}{2  L^2d}.
\end{align}
For $k \in [M-1]$ and any $\alpha \geq 0$, we have that the set 
\[
\left\{ \rho_{i}^{(t)}(\theta_k) \geq   K(\theta_M, \theta_k) + \epsilon \right\} 
\]
lies in the complement of the following set
\begin{align*}
\left\{ \rho_{i}^{(t)}(\theta_k) -\rho_{i}^{(t)}(\theta_M) <   K(\theta_M, \theta_k) + \epsilon - \alpha \right\} 
\\
\cap \left\{ \rho_{i}^{(t)}(\theta_M)<\alpha \right\},
\end{align*}
which implies
\begin{align}
&\P \left( \rho_{i}^{(t)}(\theta_k) \geq   K(\theta_M, \theta_k) + \epsilon \right)
\nonumber
\\
&\leq 
\P \left( \rho_{i}^{(t)}(\theta_k) -\rho_{i}^{(t)}(\theta_M) \geq   K(\theta_M, \theta_k) + \epsilon - \alpha\right) 
\nonumber
\\
& \hspace{0.5cm}+ \P \left( \rho_{i}^{(t)}(\theta_M) \geq\alpha \right).
\end{align}
Using Lemma~\ref{lemma:exp_k_t} we have that for every $\delta > 0$ there exists a $T$ such that for all $t \geq T$
\begin{align}
&\P \left( \rho_{i}^{(t)}(\theta_k) \geq   K(\theta_M, \theta_k) + \epsilon \right)
\nonumber
\\
&\leq \exp{\left(-\frac{ (\epsilon - \alpha)^2}{2 L^2d}t + \delta t \right)} 
\\
&+
\exp{ \left( -\min_{ k \in [M-1]} \left\{\frac{K(\theta_M, \theta_k)^2}{2 L^2d}\right\}t + \delta t \right) }.
\end{align} 
Taking limit as $\alpha \to 0^+$ for $0 < \epsilon \leq   L - K(\theta_M, \theta_k)$ we have
\begin{align}
&\lim_{t \to \infty}\frac{1}{t} \log \P 
\left(  \rho_{i}^{(t)}(\theta_k) \geq   K(\theta_M, \theta_k) + \epsilon \right)
\nonumber
\\
&\leq 
-\frac{1}{2 L^2d} \min \left\{\epsilon^2,   \min_{j \in [M-1]} K^2(\theta_M, \theta_j) \right\}.
\end{align} 
For $\epsilon \geq   L -   K(\theta_M, \theta_k)$ we have
\begin{align}
\lim_{t \to \infty}\frac{1}{t} \log 
\P \left( \rho_{i}^{(t)}(\theta_k) \geq   K(\theta_M, \theta_k) + \epsilon \right)
\nonumber
\\
\leq -\min_{k \in [M-1]} \left\{\frac{K(\theta_M, \theta_k)^2}{2 L^2d}\right\}.
\end{align}

\begin{lemma}
\label{lemma:exp_k_t}
For all $\alpha > 0$, we have the following for the sequence $\est{i}{t}(\theta_M)$
\begin{align}
\lim_{t \to \infty}\frac{1}{t} \log \P \left(  \rho_{i}^{(t)}(\theta_M) \geq \alpha \right) 
\nonumber
\\
\leq -\min_{k \in [M-1]} \left\{\frac{K(\theta_M, \theta_k)^2}{2L^2d}\right\}.
\end{align}
\end{lemma}

\begin{IEEEproof}
For any $\alpha > 0$, consider
\begin{align}
&\P \left( \rho_{i}^{(t)}(\theta_k) \geq \alpha \right)
\nonumber
\\
& \leq \sum_{k \in [M-1]} \P \left( \frac{1}{M-1} \left(1 -  e^{-\alpha t} \right) \leq \est{i}{t}(\theta_k) \right)
\nonumber
\\
& = \sum_{k \in [M-1]} \P \left( \rho_{i}^{(t)}(\theta_k) ) \leq   K(\theta_M, \theta_k) -  _t(\theta_k) \right), 
\end{align}
where $\eta_t(\theta_k) =   K(\theta_M, \theta_k) - \frac{1}{t}\log (M-1) + \frac{1}{t} \log \left(1 -  e^{-\alpha t} \right)$.
For every $\epsilon > 0$, there exists $T(\epsilon)$ such that for all $t \geq T(\epsilon)$ we have
\begin{align*}
&\P \left( \rho_{i}^{(t)}(\theta_k)  \geq \alpha\right)
\\
&\leq \sum_{k \in [M-1]} \P \left( \rho_{i}^{(t)}(\theta_k)  \leq   K(\theta_M, \theta_k) -   K(\theta_M, \theta_k) - \epsilon \right)
\\
& = \sum_{k \in [M-1]} \P \left( \rho_{i}^{(t)}(\theta_k) \leq \epsilon \right).
\end{align*}


Therefore, for every $\epsilon> 0$, $\delta > 0$, there exists $T = \max \{ T(\epsilon), T(\delta)\}$ such that for all $t \geq T$ we have
\begin{align*}
&\P \left( \rho_{i}^{(t)}(\theta_M) \geq \alpha\right)
\\
&\leq 
(M - 1) \max_{k \in [M-1]} \exp\left\{-\frac{ (  K(\theta_M, \theta_k) - \epsilon)^2}{2  L^2d}t + \delta t\right\}.
\end{align*}
By taking limit and making $\epsilon$ arbitrarily small, we have
\begin{align*}
&\lim_{t \to \infty}\frac{1}{t} \log \P \left(  \rho_{i}^{(t)}(\theta_M) \geq \alpha \right) 
\\
&\leq -\min_{k \in [M-1]} \left\{\frac{K(\theta_M, \theta_k)^2}{2L^2d}\right\}.
\end{align*}
\end{IEEEproof}

\subsubsection{Proof of Corollary~\ref{cor:concentration_learning_rate}} 
From Theorem~\ref{thm:expconvhoeffding}, we have
\begin{align*}
&\lim_{t \to \infty} \frac{1}{t} \log \P \left( \mu_i \geq   \min_{k \in [M-1]} K(\theta_M, \theta_k) + \epsilon \right)
\\
&\leq 
-\frac{1}{2  L^2d} \min \left\{\epsilon^2,   \min_{k \in [M-1]} K^2(\theta_M, \theta_k) \right\}.
\end{align*}
Now, applying Borel-Cantelli Lemma to the above equation we have
\begin{align*}
\mu_i \leq   \min_{k \in [M-1]} K(\theta_M, \theta_k) \quad \P\mbox{-a.s.}
\end{align*}
Combining this with Corollary~\ref{cor:conv} we have
\begin{align*}
\mu_i =   \min_{k \in [M-1]} K(\theta_M, \theta_k) \quad \P \mbox{-a.s.}
\end{align*} 

\subsection{Proof of Theorem~\ref{thm:ldp_vector_belief}}

\add{
To prove that $ \frac{1}{t} \log \estc{i}$ satisfies the LDP, first we establish the LDP satisfied by the following vector 
\begin{align}
\mathbf{Q}_i^{(t)}
=
\left[
\frac{\est{i}{t}(\theta_1)}{\est{i}{t}(\theta_M)},
\frac{\est{i}{t}(\theta_2)}{\est{i}{t}(\theta_M)},
\ldots,
\frac{\est{i}{t}(\theta_{M-1})}{\est{i}{t}(\theta_M)}
\right]^T.
\end{align}
Note that $\mathbf{Q}_i^{(t)} = \frac{\estc{i}}{\est{i}{t}(\theta_M)} $. From Lemma~\ref{lemma:jointLDP}, we obtain that $\frac{1}{t}\log  \mathbf{Q}_i^{(t)}$ satisfies the LDP with rate function $I(\cdot)$, as given by equation~(\ref{eq:rate_func}). Now we apply the Contraction principle [Fact~\ref{fact:contractionprinciple}], for
\begin{align*}
&\mathcal{X} = \mathbb{R}^{M-1},\quad
\mathcal{Y} = \mathbb{R}^{M-1},\\
&T\left(
\mbf{x}
\right) 
= 
g
\left(
\mbf{x}
\right),
\quad \forall \,
\mbf{x} \in \mathbb{R}^{M-1}, 
\\
&\P_t = 
\P \left( \frac{1}{t}\log  \mathbf{Q}_i^{(t)} \in \cdot \right),
\\
&\mathsf{Q}_t =
\P \left(
g \left( \frac{1}{t}\log \mathbf{Q}_i^{(t)} \right) \in \cdot
\right),
\end{align*}
and we get that $g \left( \frac{1}{t}\log \mathbf{Q}_i^{(t)} \right)$ satisfies an LDP with a rate function $J(\cdot)$, i.e., for every $F \subset \mathbb{R}^{M-1}$ we have
\begin{align}
\label{eq:fLDP_inf}
\liminf_{t \to \infty} \frac{1}{t} \log 
\P \left(
g \left( \frac{1}{t}\log \mathbf{Q}_i^{(t)} \right) \in F
\right)
\geq - \inf_{\mbf{y} \in F^o} J(\mbf{y}),
\end{align}
and
\begin{align}
\label{eq:fLDP_sup}
\limsup_{t \to \infty} \frac{1}{t} \log 
\P \left(
g \left( \frac{1}{t}\log \mathbf{Q}_i^{(t)} \right) \in F
\right)
\leq - \inf_{\mbf{y} \in \bar{F}} J(\mbf{y}).
\end{align}

Combining Lemma~\ref{lemma:f_vec_LDP} with equations~(\ref{eq:fLDP_inf}) and (\ref{eq:fLDP_sup}), we obtain that $\frac{1}{t}\log \estc{i}$ satisfies the LDP with rate function $J(\cdot)$ as well. Hence, we have the assertion of the theorem.

}

\add{
\begin{lemma}
\label{lemma:jointLDP}
The random vector $\frac{1}{t} \log \mathbf{Q}_i^{(t)}$ satisfies the LDP with rate function given by $I(\cdot)$ in \eqref{eq:rate_func_single_theta}. That is, for any set $F \subset \mathbb{R}^{M-1}$ with interior $F^o$ and closure $\bar{F}$, we have
\begin{align}
\liminf_{t \to \infty}\frac{1}{t} \log \P 
\left( 
\frac{1}{t}\log \mathbf{Q}_i^{(t)} \in  F
\right)
\geq
-\inf_{\mbf{x} \in F^o} I(\mbf{x}),
\end{align}
and
\begin{align}
\limsup_{t \to \infty}\frac{1}{t} \log \P 
\left( 
\frac{1}{t}\log \mathbf{Q}_i^{(t)} \in  F
\right)
\leq 
-\inf_{\mbf{x} \in \bar{F}} I(\mbf{x}).
\end{align}
\end{lemma}

\begin{IEEEproof}
Using the learning rule we have
\begin{align}
\frac{1}{t}\log \mathbf{Q}_i^{(t)}
&= \frac{1}{t} \sum_{\tau = 1}^{t} \sum_{j = 1}^n W^{\tau}(i,j) \mathbf{L}_j^{(t - \tau + 1)}
\nonumber
\\
&= \frac{1}{t} \sum_{\tau = 1}^{t} \sum_{j = 1}^n \left( W^{\tau}(i,j) - v_j \right) \mathbf{L}_j^{(t - \tau + 1)} 
\nonumber
\\
&\hspace{0.5cm}+ \frac{1}{t} \sum_{\tau = 1}^{t} \mbf{Y}^{(\tau)},
\end{align}
where $\mbf{L}$ is given by \eqref{eq:llr_vector} and $\mbf{Y}$ by \eqref{eq:llr_corr}.
Using Cramer's Theorem [Fact~\ref{fact:cramerstheorem}] in $\mathbb{R}^{M-1}$, for any set $F \subset \mathbb{R}^{M-1}$, we have
\begin{align}
\liminf_{t \to \infty}\frac{1}{t} \log \P 
\left( 
\frac{1}{t} \sum_{\tau = 1}^{t} \mbf{Y}^{(\tau)} \in  F
\right)
\geq
-\inf_{\mbf{x} \in F^o} I(\mbf{x}),
\end{align}
and
\begin{align}
\limsup_{t \to \infty}\frac{1}{t} \log \P 
\left( 
\frac{1}{t} \sum_{\tau = 1}^{t} \mbf{Y}^{(\tau)} \in  F
\right)
\leq
-\inf_{\mbf{x} \in \bar{F}} I(\mbf{x}).
\end{align}
Consider
\begin{align}
& \left| \frac{1}{t} \sum_{\tau = 1}^{t} \sum_{j = 1}^n \left( W^{\tau}(i,j) - v_j \right) \mathbf{L}_j^{(t - \tau + 1)} \right|
\nonumber
\\
&\leq 
\frac{  n}{t} \sum_{\tau = 1}^{t} \left| \lambda_{\max}^{\tau}(W)\right|\left( \sum_{j = 1}^n \left| \mathbf{L}_j^{(t - \tau + 1)}  \right| \right).
\end{align}
From Assumption~\ref{assume:finite_lmgf}, we have that $\Lambda(\boldsymbol{\lambda})$ is finite for $\boldsymbol{\lambda} \in \mathbb{R}^n$. Now, using Lemma~\ref{lemma:rate_error_term}, we have 
\begin{align}
\lim_{t \to \infty} \frac{1}{t}
\log \P \left(
\left| \frac{1}{t} \sum_{\tau = 1}^{t} \sum_{j = 1}^n \left( W^{\tau}(i,j) - v_j \right) \mathbf{L}_j^{(t - \tau + 1)} \right|
 \geq \boldsymbol{\delta}
\right)
\nonumber
\\
= -\infty.
\end{align}
Using Lemma~\ref{lemm:ldp_sum} on $\frac{1}{t}\log \mathbf{Q}_i^{(t)}$, we have the assertion of the theorem.
\end{IEEEproof}
}

\add{
\begin{fact}[Cramer's Theorem, Theorem 3.8~\cite{den2008large}]
\label{fact:cramerstheorem}
Consider a sequence of $d$-dimensional i.i.d random vectors $\{ \mbf{X}_n\}_{n = 1}^{\infty}$. Let $\mbf{S}_n = \frac{1}{n} \sum_{i = 1}^n \mbf{X}_i$. Then, the sequence of $\mbf{S}_n$ satisfies a large deviation principle with rate function $\Lambda^*(\cdot)$, namely: For any set $F \subset \mathbb{R}^d$,
\begin{align}
\liminf_{n \to \infty} \frac{1}{n} \log \P(\mbf{S}_n \in F)
\geq 
- \inf_{\mbf{x} \in F^o},
\end{align}
and
\begin{align}
\limsup_{n \to \infty} \frac{1}{n} \log \P(\mbf{S}_n \in F)
\leq 
- \inf_{\mbf{x} \in \bar{F}},
\end{align}
where $\Lambda^*(\cdot)$ is given by
\begin{align}
\Lambda^*(\mathbf{x}) 
\defeq
\sup_{\boldsymbol{\lambda} \in \mathbb{R}^{d}} \left\{ \left\langle \boldsymbol{\lambda}, 
\mathbf{x}
\right\rangle - \Lambda(\boldsymbol{\lambda})\right\}.
\end{align}
and $\Lambda(\cdot)$ is the log moment generating function of $\mbf{S}_n$ which is given by
\begin{align}
\Lambda(\boldsymbol{\lambda}) 
&\defeq
\log \expe [e^{\langle \boldsymbol{\lambda}, \mbf{Y}\rangle}].
\end{align}
\end{fact}
}

\add{
\begin{fact}[Contraction principle, Theorem 3.20~\cite{den2008large}]
\label{fact:contractionprinciple}
Let $\{\P_t\}$ be a sequence of probability measures on a Polish space $\mathcal{X}$ that satisfies LDP with rate function $I$. Let 
\begin{align}
\left\{
\begin{array}{ll}
	\mathcal{Y} & \text{be a Polish space}\\
	T:\mathcal{X} \to \mathcal{Y} & \text{a continuous map}\\
	\mathsf{Q}_t = \P_t \circ T^{-1} & \text{an image probability measure.}
\end{array}
\right.
\end{align}
Then $\{ \mathsf{Q}_t\}$ satisfies the LDP on $\mathcal{Y}$ with rate function $J$ given by 
\begin{align}
J(y) = \inf_{x \in \mathcal{X}: T(x) = y} I(x).
\end{align}
\end{fact}
}

\add{
\begin{lemma}
\label{lemma:f_vec_LDP}
For every set $F \subset \mathbb{R}^{M-1}$ and for all $i \in [n]$, we have
\vspace{-0.2cm}
\begin{align}
&\liminf_{t \to \infty} \frac{1}{t} \log 
\P \left( 
\frac{1}{t}\log \estc{i}
\in  F 
\right)
\nonumber
\\
&\geq \liminf_{t \to \infty} \frac{1}{t} \log 
\P \left( 
g \left( \frac{1}{t}\log \mathbf{Q}_i^{(t)} \right)  
\in  F 
\right),
\end{align}
and
\begin{align}
&\limsup_{t \to \infty} \frac{1}{t} \log 
\P \left( 
\frac{1}{t}\log \estc{i}
\in  F 
\right)
\nonumber
\\
&\leq \limsup_{t \to \infty} \frac{1}{t} \log 
\P \left( 
g \left( \frac{1}{t}\log \mathbf{Q}_i^{(t)} \right) 
\in  F 
\right).
\end{align}
\end{lemma}

\begin{IEEEproof}
For all $t \geq 0$, we have
\begin{align}
&\frac{1}{t}\log \estc{i}
=
g\left( \frac{1}{t} \log \mbf{Q}_i^{(t)} \right)
\nonumber
\\
&-\frac{1}{t} \log \left( 
e^{-C^{(t)} t} + 
\sum_{j = 1}^{M-1} e^{ g_j \left( \frac{1}{t} \log \mbf{Q}_i^{(t)} \right) t}
\right) \mbf{1},
\end{align}
where 
\begin{align*}
C^{(t)} = \max 
\left\{ 
0, 
\frac{1}{t}\log \frac{\est{i}{t}(\theta_1)}{\est{i}{t}(\theta_M)},
\frac{1}{t}\log \frac{\est{i}{t}(\theta_2)}{\est{i}{t}(\theta_M)}, \right.
\\
\left. \ldots,
\frac{1}{t}\log \frac{\est{i}{t}(\theta_{M-1})}{\est{i}{t}(\theta_M)}
\right\}.
\end{align*}
Also for all $t \geq 0 $, we have
\begin{align*}
1 \leq e^{-C^{(t)} t} + 
\sum_{j = 1}^{M-1} e^{ g_j \left( \frac{1}{t} \log \mbf{Q}_i^{(t)} \right) t} 
 \leq M.
\end{align*}
Hence for all $\epsilon > 0$, there exists $T(\epsilon)$ such that for all $t \geq T(\epsilon)$ we have
\begin{align}
g\left( \frac{1}{t} \log \mbf{Q}_i^{(t)} \right) - \epsilon \mbf{1}
 \leq
 \frac{1}{t}\log \estc{i}
\leq
g\left( \frac{1}{t} \log \mbf{Q}_i^{(t)} \right).
\end{align}
For any $F \subset \mathbb{R}^{M-1}$, let $F_{\epsilon^+} = \{\mbf{x} + \delta \mbf{1}, \forall \, 0 < \delta \leq \epsilon \text{ and }\mbf{x} \in F \}$, $F_{\epsilon^-} = \{\mbf{x} - \delta \mbf{1}, \forall \, 0 < \delta \leq \epsilon \text{ and }\mbf{x} \in F\}$. Therefore, for every $\epsilon > 0$ we have
\begin{align}
\liminf_{t \to \infty} \frac{1}{t} \log \P \left( g\left( \frac{1}{t} \log \mbf{Q}_i^{(t)} \right) \in F \right)
\nonumber
\\
\leq
\liminf_{t \to \infty} \frac{1}{t} \log \P \left( \frac{1}{t}\log \estc{i} \in F_{\epsilon^-} \right).
\end{align}
Making $\epsilon$ arbitrarily small $F_{\epsilon^-} \to F$, by monotonicity and continuity of probability measure we have
\begin{align}
\label{eq:lower_bound}
\liminf_{t \to \infty} \frac{1}{t} \log \P \left( g\left( \frac{1}{t} \log \mbf{Q}_i^{(t)} \right) \in F \right)
\nonumber
\\
\leq
\liminf_{t \to \infty} \frac{1}{t} \log \P \left( \frac{1}{t}\log \estc{i} \in F \right).
\end{align} 
For $t \geq T(\epsilon)$ we also have
\begin{align}
\frac{1}{t}\log \estc{i} 
\leq 
g\left( \frac{1}{t} \log \mbf{Q}_i^{(t)} \right)
\leq
\frac{1}{t}\log \estc{i} + \epsilon \mbf{1}.
\end{align}
This implies for every $\epsilon > 0$ we have
\begin{align}
&\limsup_{t \to \infty} \frac{1}{t} \log \P \left( \frac{1}{t}\log \estc{i} \in F \right)
\nonumber
\\
&\leq
\limsup_{t \to \infty} \frac{1}{t} \log \P \left( g\left( \frac{1}{t} \log \mbf{Q}_i^{(t)} \right) \in F_{\epsilon^+} \right).
\end{align}
Again, by making $\epsilon$ arbitrarily small we have
\begin{align}
\label{eq:upper_bound}
&\limsup_{t \to \infty} \frac{1}{t} \log \P \left( \frac{1}{t}\log \estc{i} \in F \right)
\nonumber 
\\
&\leq
\limsup_{t \to \infty} \frac{1}{t} \log \P \left( g\left( \frac{1}{t} \log \mbf{Q}_i^{(t)} \right) \in F \right).
\end{align}
Hence, we have the assertion of the lemma.
\end{IEEEproof}
}

\subsection{Proof of the Lemmas}

\add{
\begin{lemma}
\label{lemma:rate_error_term}
Let $q$ be a real number such that $q \in (0,1)$. Let $\mbf{X}_i$ be a sequence of non-negative i.i.d random vectors in $\mathbb{R}^{n}$, distributed as $\mbf{X}$ and let $\Lambda(\boldsymbol{\lambda})$ denote its log moment generating function which is finite for $\boldsymbol{\lambda} \in \mathbb{R}^n$, then for every $\delta >0$, we have 
\begin{align}
\lim_{t \to \infty} \frac{1}{t} \log \P \left(\frac{1}{t}\sum_{i =1}^t (q)^i \mbf{X}_i \geq \delta \mbf{1} \right) = -\infty.
\end{align}
\end{lemma}


\begin{IEEEproof}
Applying Chebycheff's inequality and using the definition of log moment generating function, for $\boldsymbol{\lambda} \in \mathbb{R}^n$, we have
\begin{align}
\P \left(\frac{1}{t}\sum_{i =1}^t (q)^i \mbf{X}_i \geq \delta \mbf{1} \right) 
&\leq
e^{ - t\left( \langle \boldsymbol{\lambda}, \delta \mbf{1}\rangle - \frac{1}{t} \sum_{i = 1}^{t} \Lambda((q)^i \boldsymbol{\lambda}) \right)}. 
\end{align}
From convexity of $\Lambda$, we have 
$
\sum_{i = 1}^t  \Lambda ((q)^i \boldsymbol{\lambda}) 
\leq 
\Lambda(\boldsymbol{\lambda}) \sum_{i = 1}^t (q)^i.
$
Since $\Lambda(\boldsymbol{\lambda})$ is finite and $\sum_{i = 1}^{\infty} (q)^i < \infty$, for all $\delta > 0$ we have
\begin{align}
\lim_{t \to \infty} \frac{1}{t} \log
\P\left(\frac{1}{t}\sum_{i =1}^t (q)^i \mbf{X}_i \geq \delta \mbf{1} \right)
\leq
-  \langle \boldsymbol{\lambda}, \delta \mbf{1}\rangle.
\end{align}
Since, the above equation is true for all $\boldsymbol{\lambda} \in \mathbb{R}^n$, we have the assertion of the lemma.
\end{IEEEproof}
}

\add{
\begin{lemma}
\label{lemm:ldp_sum}
Consider a sequence $\{\mbf{Z}^{(t)} \}_{t = 0}^{\infty}$ where $\mbf{Z}^{(t)} \in \mathbb{R}^d$ such that
\begin{align}
\mbf{Z}^{(t)} = \mbf{X}^{(t)} + \mbf{Y}^{(t)},
\end{align}
where sequences $\{\mbf{X}^{(t)} \}_{t = 0}^{\infty}$ and $\{\mbf{Y}^{(t)}\}_{t = 0}^{\infty}$ have the following properties:
\begin{enumerate}

\item
The sequence $\{\mbf{X}^{(t)}\}_{t = 0}^{\infty}$ satisfies
\begin{align}
&\liminf_{t \to \infty} \frac{1}{t} \log \P \left(\mbf{X}^{(t)} \in F \right) \geq - \inf_{\mbf{x} \in F^{o}}I_{X}(\mbf{x}), \\
&\limsup_{t \to \infty} \frac{1}{t} \log \P \left(\mbf{X}^{(t)} \in F \right) \leq - \inf_{\mbf{x} \in \bar{F}}I_{X}(\mbf{x}),
\end{align}
where $I_{X}: \mathbb{R}^d \to \mathbb{R}$ is a well-defined LDP rate function.

\item
For every $\epsilon > 0$, sequence $\{\mbf{Y}^{(t)} \}_{t = 0}^{\infty}$ satisfies
\begin{align}
&\lim_{t \to \infty} \frac{1}{t} \log \P (|\mbf{Y}^{(t)}| \geq \epsilon \mbf{1}) = -\infty.
\end{align}

\end{enumerate} 

Then $\{\mbf{Z}^{(t)} \}_{t = 0}^{\infty}$ satisfies
\begin{align}
&\liminf_{t \to \infty} \frac{1}{t} \log \P (\mbf{Z}^{(t)} \in F) \geq - \inf_{\mbf{x} \in F^{o}}I_{X}(\mbf{x}), \\
&\limsup_{t \to \infty} \frac{1}{t} \log \P (\mbf{Z}^{(t)} \in F)\leq - \inf_{\mbf{x} \in \bar{F}}I_{X}(\mbf{x}).
\end{align}
\end{lemma}
}

\add{
\begin{IEEEproof}
For every $t \geq 0$, we have
\begin{align*}
&\P \left( \mbf{Z}^{(t)} \in F_{\epsilon^+} \cup F_{\epsilon^-} \right)
\\
&\geq
\P\left( \{ \mbf{X}^{(t)} \in F \} \cap \{ |\mbf{Y}^{(t)}| \leq \epsilon \mbf{1}\} \right)
\\
&\geq
\P\left(  \mbf{X}^{(t)} \in F \right)
-
\P \left(  |\mbf{Y}^{(t)}| > \epsilon \mbf{1} \right).
\end{align*}
For all $\delta > 0$, there exists a $T(\delta)$ such that for all $t \geq T(\delta)$ we have
\begin{align*}
\P\left(  \mbf{X}^{(t)} \in F \right)
\geq e^{-\inf_{\mbf{x} \in F^{o}}I_{X}(\mbf{x})t - \delta t}.
\end{align*}
For all $B > 0$,  there exists a $T(B)$ such that for all $t \geq T(B)$ we have
\begin{align*}
\P \left(  |\mbf{Y}^{(t)}| > \epsilon \mbf{1} \right) \geq e^{-B t}.
\end{align*}
Now choose $B > \inf_{\mbf{x} \in F^{o}}I_{X}(\mbf{x})+ \delta$ and $t \geq \max \{T(\delta), T(B) \}$, then we have
\begin{align*}
&\P \left( \mbf{Z}^{(t)} \in F_{\epsilon^+} \cup F_{\epsilon^-} \right)
\\
&\geq
e^{-\inf_{\mbf{x} \in F^{o}}I_{X}(\mbf{x})t - \delta t} \left( 1- e^{ -Bt + \inf_{\mbf{x} \in F^{o}}I_{X}(\mbf{x})t + \delta t} \right).
\end{align*}
Sending $\epsilon$ to zero and taking the limit we have
\begin{align*}
\liminf_{t \to \infty} \frac{1}{t}\log \P \left( \mbf{Z}^{(t)} \in F \right)
\geq
-\inf_{\mbf{x} \in F^{o}}I_{X}(\mbf{x}).
\end{align*}
Similarly, using the fact that $\P(\{\mbf{Z}^{(t)} \in F \} \cap \{|\mbf{Y}^{(t)}| \leq \epsilon \mbf{1} \} ) \leq \P \left( \mbf{X}^{(t)} \in F_{\epsilon^+}\right)$ we have the other LDP bound.
\end{IEEEproof}
}

\bibliographystyle{IEEEtran}
\bibliography{social_learning}

\end{document}